# Robust reliable control of uncertain 2D discrete switched systems with state delays

**Shipei Huang,** **Zhengrong Xiang***

School of Automation, Nanjing University of Science and Technology

Nanjing, 210094, People's Republic of China

*corresponding author, e-mail: xiangzr@mail.njust.edu.cn

**Abstract:** This paper is concerned with the problem of robust reliable control for a class of uncertain 2D discrete switched systems with state delays represented by a model of Roesser type. The parameter uncertainties are assumed to be norm-bounded. Firstly, delay-dependent sufficient condition for the exponential stability of the discrete 2D systems with state delays is established. Then, the concept of average dwell time is extended to 2D switched systems, and a reliable state feedback controller is developed in terms of linear matrix inequalities (LMIs) such that the resulting closed-loop system is exponentially stable for all admissible uncertainties and actuator failures. The dwell time approach is utilized for the stability analysis and controller design. Finally, an example is included to demonstrate the effectiveness of the proposed approach.

**Keywords:** 2D discrete systems; switched systems; state delays; reliable control; dwell time

## 1. Introduction

Two-dimensional (2D) systems have received considerable attention over the past few decades due to their wide applications in many areas such as multi-dimensional digital filtering, linear image processing, signal processing, and process control [1-3]. The stability analysis of 2-D discrete systems has attracted a great deal of interest and some significant results have been obtained [4-12]. Furthermore, time delays frequently occur in practical systems and are often the source of instability, there are many examples containing inherent delays in practical 2-D discrete systems, the stability of 2-D discrete systems with state delays have also been studied in [8-12].

On the other hand, switched systems have also been attracted considerable attention during the past several decades [13-23]. A switched system is a hybrid system which consists of a finite number of continuous-time or discrete-time subsystems and a switching signal specifying the switch between these subsystems. This class of systems has numerous applications in many fields, such as mechanical systems, the automotive industry, aircraft and air traffic control, switched power converters. So far, dwell time approach is an important and effective approach to studying switched systems. Recently, the dwell time approach is applied widely to deal with switched systems, see, for example, [14], [15], [22] and references therein.

However, the switch phenomenon may also occur in 2D discrete systems, and study on these 2D discrete switched systems will also be significant. There are a few reports on 2D discrete switched systems at present. Benzaouia *et al.* [24] firstly consider 2D switched systems with arbitrary switched sequences, and the process of switched is considered as a Markovian jumping one. Furthermore, the stabilizability problem of 2D discrete switched systems is investigated in [25], a sufficient condition for the asymptotic stability of such systems is proposed and a stabilizing controller is developed in terms of LMI. It should be noted that these papers focus on studying the asymptotical stability of the 2D switched systems, and the obtained results are based on common and multiple Lyapunov function approachs, the problem of stability for 2D switched systems via the dwell time approach has not been investigated to date, especially for the exponential stability problem of 2D switched systems with state delay, which motivates us to shorten such a gap in the present investigation.

It is well known that the actuators may be subjected to failures in practical operation. Therefore, it is of practical interest to design a control system which can tolerate faults of actuators. A control system is said to be reliable if it retains certain properties when there exist failures. It should be noted that in normal cases, a controller with fixed gain is easily implemented, and could meet the requirement in practical applications. But when failure occurs, the conventional controller will become conservative and may not satisfy certain control performance indexes. Reliable control is a kind of effective control approach to improve system reliability. Several approaches for designing reliable controllers have been proposed, and some of which have been used to research the problem of reliable control for switched systems [26-30]. To the best of our knowledge, there are few works about designing robust reliable controller for 2D discrete switched systems with state delays, which motivates our present study.

In this paper, we are interested in designing a reliable stabilizing controller for 2D discrete switched systems with state delays such that the closed-loop system is exponentially stable. The dwell time approach is utilized for the stability analysis and controller design. The remainder of the paper is organized as follows. In Section 2, problem formulation and some necessary lemmas are given. In Section 3, based on the dwell time approach, stability for 2D discrete switched systems with state delays are addressed, and a delay-dependent sufficient condition for the existence of a robust reliable controller is derived in terms of a set of matrix inequalities. A numerical example is provided to illustrate the effectiveness of the proposed approach in Section 4. Concluding remarks are given in Section 5.

**Notations**: Throughout this paper, the superscript "$T$" denotes the transpose, and the notation $X \geq Y\left(X > Y\right)$ means that matrix $X - Y$ is positive semi-definite (positive definite, respectively). $\|*\|$ denotes the Euclidean norm. $I$ represents identity matrix with appropriate dimension. $I_h$ is the identity matrix with $n_1$ appropriate dimension and $I_v$ is the identity matrix with $n_2$ appropriate dimension. $diag\{a_i\}$ denotes diagonal matrix with the diagonal elements $a_i, i = 1, 2, ..., n$. $X^{-1}$ denotes the inverse of $X$. The asterisk $*$ in a matrix is used to denote term that is induced by symmetry. The set of all nonnegative integers is represented by $Z_+$.

## 2. Problem formulation and preliminaries

Consider the following uncertain 2D discrete linear switched systems with state delays described by the Roesser model:

$$\begin{bmatrix} x^h\left(i+1, j\right) \\ x^v\left(i, j+1\right) \end{bmatrix} = \hat{A}^{\sigma(m)} \begin{bmatrix} x^h\left(i, j\right) \\ x^v\left(i, j\right) \end{bmatrix} + \hat{A}_d^{\sigma(m)} \begin{bmatrix} x^h\left(i - d_h\left(i\right), j\right) \\ x^v\left(i, j - d_v\left(j\right)\right) \end{bmatrix} + B^{\sigma(m)} u^f\left(i, j\right) \tag{1}$$

where $x^h\left(i, j\right)$ is the horizontal state in $R^{n_1}$, $x^v\left(i, j\right)$ is the vertical state in $R^{n_2}$, $x\left(i, j\right)$ is the whole state in $R^n$ with $n = n_1 + n_2$. $u^f\left(i, j\right)$ is the control input of actuator fault, $m=i+j$ is the sum of $i$ and $j$. $\sigma(m)$ is a switching rule which takes its values in the finite set $\underline{N}:=\{1,\cdots,N\}$, and $\sigma(m)=k$ means that the $k$th subsystem is activated. $N$ is the number of subsystems, $i$ and $j$ are integers in $Z_+$. $d_h(i)$ and $d_v(j)$ are delays along horizontal and vertical directions, respectively. We assume $d_h(i)$ and $d_v(j)$ satisfying

$$d_{hL} \le d_h(i) \le d_{hH},\ d_{vL} \le d_v(j) \le d_{vH} \tag{2}$$

where $d_{hL}$, $d_{hH}$ and $d_{vL}$, $d_{vH}$ denote the lower and upper delay bounds along horizontal and vertical directions, respectively.

$\hat{A}^{\sigma(m)}$, $\hat{A}_d^{\sigma(m)}$ are uncertain real-valued matrices with appropriate dimensions which are assumed to be of the form

$$\hat{A}^{\sigma(m)} = A^{\sigma(m)} + H^{\sigma(m)} F^{\sigma(m)}(i,j) E^{\sigma(m)}$$

$$\hat{A}_d^{\sigma(m)} = A_d^{\sigma(m)} + H^{\sigma(m)} F^{\sigma(m)}(i,j) E_d^{\sigma(m)} \tag{3}$$

with

$$A^{\sigma(m)} = \begin{bmatrix} A_{11}^{\sigma(m)} & A_{12}^{\sigma(m)} \\ A_{21}^{\sigma(m)} & A_{22}^{\sigma(m)} \end{bmatrix},\ A_d^{\sigma(m)} = \begin{bmatrix} A_{d11}^{\sigma(m)} & A_{d12}^{\sigma(m)} \\ A_{d21}^{\sigma(m)} & A_{d22}^{\sigma(m)} \end{bmatrix},$$

$$H^{\sigma(m)} = \begin{bmatrix} H_1^{\sigma(m)} \\ H_2^{\sigma(m)} \end{bmatrix},\ E^{\sigma(m)} = \begin{bmatrix} E_1^{\sigma(m)} \\ E_2^{\sigma(m)} \end{bmatrix},\ E_d^{\sigma(m)} = \begin{bmatrix} E_{d1}^{\sigma(m)} \\ E_{d2}^{\sigma(m)} \end{bmatrix},$$

where matrices $A_{11}^{\sigma(m)} \in R^{n_1\times n_1}$, $A_{12}^{\sigma(m)} \in R^{n_1\times n_2}$, $A_{21}^{\sigma(m)} \in R^{n_2\times n_1}$, $A_{22}^{\sigma(m)} \in R^{n_2\times n_2}$, $A_{d11}^{\sigma(m)} \in R^{n_1\times n_1}$, $A_{d12}^{\sigma(m)} \in R^{n_1\times n_2}$, $A_{d21}^{\sigma(m)} \in R^{n_2\times n_1}$, $A_{d22}^{\sigma(m)} \in R^{n_2\times n_2}$, $H_1^{\sigma(m)} \in R^{n_1\times p}$, $H_2^{\sigma(m)} \in R^{n_2\times p}$, $E_1^{\sigma(m)} \in R^{q\times n_1}$, $E_2^{\sigma(m)} \in R^{q\times n_2}$ are constant matrices. $F^{\sigma(m)}(i,j) \in R^{p\times q}$ is an unknown matrix representing parameter uncertainty which satisfies

$$F^{\sigma(m)T}(i,j) F^{\sigma(m)}(i,j) \le I \tag{4}$$

The control input of actuator fault $u^f(i,j)$ can be described as

$$u^f(i,j) = \Omega^{\sigma(m)} u(i,j) \tag{5}$$

where $u(i,j) = K^{\sigma(m)} x(i,j)$ is the control input to be designed, $\Omega^k\ (k \in \underline{N})$ are the actuator fault matrices with the following form

$$\Omega^k = diag\left\{\omega_{k1}, \omega_{k2}, ..., \omega_{kl}, ..., \omega_{kn}\right\} \tag{6}$$

where $0 \le \omega_{Lkl} \le \omega_{kl} \le \omega_{Hkl}$.

For simplicity, we define

$$\Omega_0^k = diag\left\{\tilde{\omega}_{k1}, \tilde{\omega}_{k2}, ..., \tilde{\omega}_{kl}, ..., \tilde{\omega}_{kn}\right\},\ \tilde{\omega}_{kl} = \frac{1}{2}\left(\omega_{Lkl} + \omega_{Hkl}\right) \tag{7}$$

$$\Xi^k = diag\left\{\xi_{k1}, \xi_{k2}, ..., \xi_{kl}, ..., \xi_{kn}\right\},\ \xi_{kl} = \frac{\omega_{Hkl} - \omega_{Lkl}}{\omega_{Hkl} + \omega_{Lkl}} \tag{8}$$

$$\Theta^k = diag\left\{\Theta_{k1}, \Theta_{k2}, ..., \Theta_{kl}, ..., \Theta_{kn}\right\},\ \Theta_{kl} = \frac{\omega_{kl} - \tilde{\omega}_{kl}}{\tilde{\omega}_{kl}} \tag{9}$$

Thus, we have

$$\Omega^k = \Omega_0^k\left(I + \Theta^k\right),\ \left|\Theta^k\right| \le \Xi^k \le I \tag{10}$$

where $\left|\Theta^k\right| = diag\left\{\left|\Theta_{k1}\right|, \left|\Theta_{k2}\right|, ..., \left|\Theta_{kl}\right|, ..., \left|\Theta_{kn}\right|\right\}$.

The boundary conditions are defined by

$$\begin{aligned}
&x^h(i,j) = h_{ij}, \quad \forall 0 \le j \le z_1,\ -d_{hH} \le i \le 0 \\
&x^h(i,j) = 0, \quad \forall j > z_1, \quad -d_{hH} \le i \le 0 \\
&x^v(i,j) = v_{ij}, \quad \forall 0 \le i \le z_2,\ -d_{vH} \le j \le 0 \\
&x^v(i,j) = 0, \quad \forall i > z_2, \quad -d_{vH} \le j \le 0
\end{aligned} \tag{11}$$

where $z_1 < \infty$ and $z_2 < \infty$ are positive integers, $h_{ij}$ and $v_{ij}$ are given vectors.

**Remark 1** In the paper, it is assumed that switch occurs only at each sampling of $i$ or $j$, here is $m = i + j$. The switch sequences of the system can be described as

$$\left(m_1, \sigma\left(m_1\right)\right), \left(m_2, \sigma\left(m_2\right)\right), ..., \left(m_\kappa, \sigma\left(m_\kappa\right)\right)... \tag{12}$$

with $\kappa \in Z_+$ $m_\kappa = i^\kappa + j^\kappa$, $m_\kappa$ denotes the $\kappa - th$ switching instant and $\sigma\left(m_\kappa\right) \in \underline{N}$ denotes the value of switch signal.

**Definition 1** The system (1) is said to be exponentially stable under $\sigma(m)$ if for a given $z \ge 0$, there exists a positive constant $c$, such that

$$\sum_{i+j=D} \left\| x(i,j) \right\|^2 \le \eta e^{-c(D-z)} \sum_{i+j=z} \left\| x(i,j) \right\|_C^2 \tag{13}$$

holds for all $D \ge z$ and a positive constants $\eta$ and

$$\sum_{i+j=z} \left\| x(i,j) \right\|_C^2 \triangleq \sup_{\substack{-d_{hH} < \theta_h \le 0, \\ -d_{vH} < \theta_v \le 0}} \sum_{i+j=z} \Big\{ \left\| x^h (i-\theta_h, j) \right\|^2 + \left\| x^v (i, j-\theta_v) \right\|^2 ,$$

$$\left\| \eta^h (i-\theta_h, j) \right\|^2 + \left\| \eta^v (i, j-\theta_v) \right\|^2 \Big\}$$

where

$$\eta^h (i-\theta_h, j) = x^h (i-\theta_h + 1, j) - x^h (i-\theta_h, j)$$

$$\eta^v (i, j-\theta_v) = x^v (i, j-\theta_v + 1) - x^v (i, j-\theta_v)$$

**Remark 2** From the definition 1, it is easy to see that when $z$ is given, $\sum_{i+j=z} \left\| x(i,j) \right\|_C^2$ will be bounded and as $D$ goes to infinity, $\sum_{i+j=D} \left\| x(i,j) \right\|^2$ will tend to zero exponentially, which also means $\left\| x(i,j) \right\|$ tends to zero.

**Definition 2** For any $i + j = D \ge z = i_z + j_z$, let $N_{\sigma(m)}(z, D)$ denote the switch number of $\sigma(m)$ on an interval $(z, D)$. If

$$N_{\sigma(m)}(z, D) \le N_0 + \frac{D-z}{\tau_a} \tag{14}$$

hold for given $N_0 \ge 0$, $\tau_a \ge 0$, then the constant $\tau_a$ is called the average dwell time and $N_0$ is the chatter bound.

**Lemma 1** [31] For a given matrix $S = \begin{bmatrix} S_{11} & S_{12} \\ S_{12}^T & S_{22} \end{bmatrix}$, where $S_{11}$, $S_{22}$ are square matrices, then the following conditions are equality.

**(i)** $S < 0$

**(ii)** $S_{11} < 0$, $S_{22} - S_{12}^T S_{11}^{-1} S_{12} < 0$

**(iii)** $S_{22} < 0$, $S_{11} - S_{12} S_{22}^{-1} S_{12}^T < 0$

**Lemma 2** [32] Let $U$, $V$, $W$ and $X$ be real matrices of appropriate dimensions with $X$ satisfying $X = X^T$, then for all $V^T V \le I$, $X + UVW + W^T V^T U^T < 0$, if and only if there exists a scalar $\varepsilon$ such that $X+\varepsilon UU^T+\varepsilon^{-1}W^TW<0$.

**Lemma 3** [33] For matrices $R_1$, $R_2$ with appropriate dimensions, there exists a positive scalar $\varepsilon$, such that

$$R_1\Sigma R_2+R_2^T\Sigma^T R_1^T\le\varepsilon R_1UR_1^T+\varepsilon^{-1}R_2^TUR_2$$

holds, where $\Sigma$ is a diagonal matrix and $U$ is known real-value matrix satisfying $|\Sigma|\le U$.

# 3. Main results

## 3.1 Stability analysis

In this subsection, we first investigate the problem of stability analysis for the 2D discrete linear non-switched systems with state delays.

**Lemma 4** Consider the following 2D discrete system with state delays

$$\begin{bmatrix}x^h(i+1,j)\\x^v(i,j+1)\end{bmatrix}=A\begin{bmatrix}x^h(i,j)\\x^v(i,j)\end{bmatrix}+A_d\begin{bmatrix}x^h(i-d_h(i),j)\\x^v(i,j-d_v(j))\end{bmatrix}\tag{15}$$

where $A$ and $A_d$ are constant matrices with appropriate dimensions. The boundary conditions are given in (11). For a given positive scalar $\alpha$, if there exist positive definite symmetric matrices $P=diag\{P_h,P_v\}$, $Q_1=diag\{Q_{1h},Q_{1v}\}$, $Q_2=diag\{Q_{2h},Q_{2v}\}$, $W_1=diag\{W_{1h},W_{1v}\}$, $W_2=diag\{W_{2h},W_{2v}\}$, $X_h=\begin{bmatrix}X_{11h}&X_{12h}\\ *&X_{22h}\end{bmatrix}$, $X_v=\begin{bmatrix}X_{11v}&X_{12v}\\ *&X_{22v}\end{bmatrix}$, $Y_h=\begin{bmatrix}Y_{11h}&Y_{12h}\\ *&Y_{22h}\end{bmatrix}$, $Y_v=\begin{bmatrix}Y_{11v}&Y_{12v}\\ *&Y_{22v}\end{bmatrix}$, and any matrices $M_h=\begin{bmatrix}M_1^h\\M_2^h\end{bmatrix}$, $N_h=\begin{bmatrix}N_1^h\\N_2^h\end{bmatrix}$, $S_h=\begin{bmatrix}S_1^h\\S_2^h\end{bmatrix}$, $M_v=\begin{bmatrix}M_1^v\\M_2^v\end{bmatrix}$, $N_v=\begin{bmatrix}N_1^v\\N_2^v\end{bmatrix}$ and $S_v=\begin{bmatrix}S_1^v\\S_2^v\end{bmatrix}$ with appropriate dimensions, such that

$$\Phi=\begin{bmatrix}\Phi_1&\Phi_2^T&\Phi_2^T&\Phi_3^TP\\ *&-(\Lambda_1)^{-1}W_1^{-1}&0&0\\ *&*&-(\Lambda_2)^{-1}W_2^{-1}&0\\ *&*&*&-P\end{bmatrix}<0\tag{16}$$

$$\Psi=diag\{\Psi_{1h},\Psi_{2h},\Psi_{3h},\Psi_{1v},\Psi_{2v},\Psi_{3v}\}\ge 0\tag{17}$$

where

$$\Phi_1 = \begin{bmatrix} \Phi_{11} & \Phi_{12} & \Phi_{13} & \Phi_{14} \\ * & \Phi_{22} & \Phi_{23} & \Phi_{24} \\ * & * & \Phi_{33} & 0 \\ * & * & * & \Phi_{44} \end{bmatrix}$$

$\Lambda_1 = diag\left\{d_{hL}I_h, d_{vL}I_v\right\}$, $\Lambda_2 = diag\left\{\left(d_{hH} - d_{hL}\right)I_h, \left(d_{vH} - d_{vL}\right)I_v\right\}$,

$\Lambda_3 = diag\left\{\alpha^{d_{hL}} I_h, \alpha^{d_{vL}} I_v\right\}$, $\Lambda_4 = diag\left\{\alpha^{d_{hH}} I_h, \alpha^{d_{vH}} I_v\right\}$,

$\Phi_{11} = Q_1 + Q_2 + \Lambda_3\left(M_1 + M_1^T + \Lambda_1 X_{11}\right) + \Lambda_2\Lambda_4 Y_{11} - \alpha P$, $\Phi_{23} = -\Lambda_3 M_2 + \Lambda_4 S_2$,

$\Phi_{12} = \Lambda_3\left(\Lambda_1 X_{12} + M_2^T\right) + \Lambda_4\left(N_1 - S_1 + \Lambda_2 Y_{12}\right)$, $\Phi_{13} = -\Lambda_3 M_1 + \Lambda_4 S_1$, $\Phi_{14} = -\Lambda_4 N_1$

$\Phi_{22} = \Lambda_1\Lambda_3 X_{22} + \Lambda_4\left(N_2 + N_2^T - S_2 - S_2^T + \Lambda_2 Y_{22}\right)$, $\Phi_{24} = -\Lambda_4 N_2$, $\Phi_{33} = -\Lambda_3 Q_1$,

$\Phi_{44} = -\Lambda_4 Q_2$, $\Phi_2 = \begin{bmatrix} A - I & A_d & 0 & 0 \end{bmatrix}$, $\Phi_3 = \begin{bmatrix} A & A_d & 0 & 0 \end{bmatrix}$,

$X_{11} = diag\left\{X_{11}^h, X_{11}^v\right\}$, $X_{12} = diag\left\{X_{12}^h, X_{12}^v\right\}$, $X_{22} = diag\left\{X_{22}^h, X_{22}^v\right\}$,

$Y_{11} = diag\left\{Y_{11}^h, Y_{11}^v\right\}$, $Y_{12} = diag\left\{Y_{12}^h, Y_{12}^v\right\}$, $Y_{22} = diag\left\{Y_{22}^h, Y_{22}^v\right\}$,

$M_1 = diag\left\{M_1^h, M_1^v\right\}$, $M_2 = diag\left\{M_2^h, M_2^v\right\}$, $N_1 = diag\left\{N_1^h, N_1^v\right\}$,

$N_2 = diag\left\{N_2^h, N_2^v\right\}$, $S_1 = diag\left\{S_1^h, S_1^v\right\}$, $S_2 = diag\left\{S_2^h, S_2^v\right\}$,

$\Psi_{1h} = \begin{bmatrix} X_h & M_h \\ * & W_{1h} \end{bmatrix}$, $\Psi_{2h} = \begin{bmatrix} Y_h & N_h \\ * & W_{2h} \end{bmatrix}$, $\Psi_{3h} = \begin{bmatrix} Y_h & S_h \\ * & W_{2h} \end{bmatrix}$,

$\Psi_{1v} = \begin{bmatrix} X_v & M_v \\ * & W_{1v} \end{bmatrix}$, $\Psi_{2v} = \begin{bmatrix} Y_v & N_v \\ * & W_{2v} \end{bmatrix}$, $\Psi_{3v} = \begin{bmatrix} Y_v & S_v \\ * & W_{2v} \end{bmatrix}$.

Then, along the trajectory of systems (15), there holds the following inequality

$$\sum_{i+j=D} V(i,j) < \alpha^{D-z} \sum_{i+j=z} V(i,j) \tag{18}$$

**Proof** Consider Lyapunov-Krasovskii functional candidate

$$V\left(x(i,j)\right) = V^h\left(x^h(i,j)\right) + V^v\left(x^v(i,j)\right) \tag{19}$$

where

$$V^h\left(x^h(i,j)\right) = \sum_{g=1}^{3} V_g^h\left(x^h(i,j)\right)$$

$$V_1^h\left(x^h(i,j)\right) = x^h(i,j)^T P_h x^h(i,j)$$

$$V_2^h\left(x^h(i,j)\right)=\sum_{r=i-d_{hL}}^{i-1} x^h(r,j)^T Q_{1h} x^h(r,j)\alpha^{i-r-1}$$

$$+\sum_{r=i-d_{hH}}^{i-1} x^h(r,j)^T Q_{2h} x^h(r,j)\alpha^{i-r-1}$$

$$V_3^h\left(x^h(i,j)\right)=\sum_{s=-d_{hL}}^{-1}\sum_{r=i+s}^{i-1}\eta^h(r,j)^T W_{1h}\eta^h(r,j)\alpha^{i-r-1}$$

$$+\sum_{s=-d_{hH}}^{-d_{hL}-1}\sum_{r=i+s}^{i-1}\eta^h(r,j)^T W_{2h}\eta^h(r,j)\alpha^{i-r-1}$$

$$V^v\left(x^v(i,j)\right)=\sum_{g=1}^{3}V_g^v\left(x^v(i,j)\right)$$

$$V_1^v\left(x^v(i,j)\right)=x^v(i,j)^T P_v x^v(i,j)$$

$$V_2^v\left(x^v(i,j)\right)=\sum_{t=j-d_{vL}}^{j-1} x^v(i,t)^T Q_{1v} x^v(i,t)\alpha^{j-t-1}$$

$$+\sum_{t=j-d_{vH}}^{j-1} x^v(i,t)^T Q_{2v} x^v(i,t)\alpha^{j-t-1}$$

$$V_3^v\left(x^v(i,j)\right)=\sum_{s=-d_{vL}}^{-1}\sum_{t=j+s}^{j-1}\eta^v(i,t)^T W_{1v}\eta^v(i,t)\alpha^{j-t-1}$$

$$+\sum_{s=-d_{vH}}^{-d_{vL}-1}\sum_{t=j+s}^{j-1}\eta^v(i,t)^T W_{2v}\eta^v(i,t)\alpha^{j-t-1}$$

$$\eta^h(r,j)=x^h(r+1,j)-x^h(r,j)$$

$$\eta^v(i,t)=x^v(i,t+1)-x^v(i,t)$$

Then we have

$$V^h\left(x^h(i+1,j)\right)-\alpha V^h\left(x^h(i,j)\right)+V^v\left(x^v(i,j+1)\right)-\alpha V^v\left(x^v(i,j)\right)$$

$$=\sum_{g=1}^{3}\left[V_g^h\left(x^h(i+1,j)\right)-\alpha V_g^h\left(x^h(i,j)\right)\right]+\sum_{g=1}^{3}\left[V_g^v\left(x^v(i,j+1)\right)-\alpha V_g^v\left(x^v(i,j)\right)\right] \quad (20)$$

with

$$V_1^h\left(x^h(i+1,j)\right)-\alpha V_1^h\left(x^h(i,j)\right)$$

$$=x^h(i+1,j)^T P_h x^h(i+1,j)-\alpha x^h(i,j)^T P_h x^h(i,j)$$

$$
\begin{aligned}
&V_2^h\left(x^h\left(i+1, j\right)\right)-\alpha V_2^h\left(x^h\left(i, j\right)\right)\\
&=\sum_{r=i-d_{hL}}^{i} x^h\left(r, j\right)^T Q_{1h} x^h\left(r, j\right)\alpha^{i-r}+\sum_{r=i-d_{hH}}^{i} x^h\left(r, j\right)^T Q_{2h} x^h\left(r, j\right)\alpha^{i-r}\\
&-\alpha\sum_{r=i-d_{hL}}^{i-1} x^h\left(r, j\right)^T Q_{1h} x^h\left(r, j\right)\alpha^{i-r-1}-\alpha\sum_{r=i-d_{hH}}^{i-1} x^h\left(r, j\right)^T Q_{2h} x^h\left(r, j\right)\alpha^{i-r-1}\\
&=x^h\left(i, j\right)^T\left(Q_{1h}+Q_{2h}\right)x^h\left(i, j\right)-\alpha^{d_{hL}} x^h\left(i-d_{hL}, j\right)^T Q_{1h} x^h\left(i-d_{hL}, j\right)\\
&-\alpha^{d_{hH}} x^h\left(i-d_{hH}, j\right)^T Q_{2h} x^h\left(i-d_{hH}, j\right)
\end{aligned}
$$

$$
\begin{aligned}
&V_3^h\left(x^h\left(i+1, j\right)\right)-\alpha V_3^h\left(x^h\left(i, j\right)\right)\\
&\le d_{hL}\eta^h\left(i, j\right)^T W_{1h}\eta^h\left(i, j\right)\\
&+\left(d_{hH}-d_{hL}\right)\eta^h\left(i, j\right)^T W_{2h}\eta^h\left(i, j\right)\\
&-\alpha^{d_{hL}}\sum_{r=i-d_{hL}}^{i-1}\eta^h\left(r, j\right)^T W_{1h}\eta^h\left(r, j\right)\\
&-\alpha^{d_{hH}}\sum_{r=i-d_{hH}}^{i-d_h(i)-1}\eta^h\left(r, j\right)^T W_{2h}\eta^h\left(r, j\right)\\
&-\alpha^{d_{hH}}\sum_{r=i-d_h(i)}^{i-d_{hL}-1}\eta^h\left(r, j\right)^T W_{2h}\eta^h\left(r, j\right)
\end{aligned}
$$

$$
\begin{aligned}
&V_1^v\left(x^v\left(i, j+1\right)\right)-\alpha V_1^v\left(x^v\left(i, j\right)\right)\\
&=x^v\left(i, j+1\right)^T P_v x^v\left(i, j+1\right)-\alpha x^v\left(i, j\right)^T P_v x^v\left(i, j\right)
\end{aligned}
$$

$$
\begin{aligned}
&V_2^v\left(x^v\left(i, j+1\right)\right)-\alpha V_2^v\left(x^v\left(i, j\right)\right)\\
&=\sum_{t=j-d_{vL}}^{j} x^v\left(i,t\right)^T Q_{1v} x^v\left(i,t\right)\alpha^{j-t}+\sum_{t=j-d_{vH}}^{j} x^v\left(i,t\right)^T Q_{2v} x^v\left(i,t\right)\alpha^{j-t}\\
&-\alpha\sum_{t=j-d_{vL}}^{j-1} x^v\left(i,t\right)^T Q_{1v} x^v\left(i,t\right)\alpha^{j-t-1}-\alpha\sum_{t=j-d_{vH}}^{j-1} x^v\left(i,t\right)^T Q_{2v} x^v\left(i,t\right)\alpha^{j-t-1}\\
&=x^v\left(i, j\right)^T\left(Q_{1v}+Q_{2v}\right)x^v\left(i, j\right)-\alpha^{d_{vL}} x^v\left(i, j-d_{vL}\right)^T Q_{1v} x^v\left(i, j-d_{vL}\right)\\
&-\alpha^{d_{vH}} x^v\left(i, j-d_{vH}\right)^T Q_{2v} x^v\left(i, j-d_{vH}\right)
\end{aligned}
$$

$$\begin{aligned}
&V_3^v\left(x^v(i,j+1)\right)-\alpha V_3^v\left(x^v(i,j)\right)\\
&\le d_{vL}\eta^v(i,j)^T W_{1v}\eta^v(i,j)\\
&+\left(d_{vH}-d_{vL}\right)\eta^v(i,j)^T W_{2v}\eta^v(i,j)\\
&-\alpha^{d_{vL}}\sum_{t=j-d_{vL}}^{j-1}\eta^v(i,t)^T W_{1v}\eta^v(i,t)\\
&-\alpha^{d_{vH}}\sum_{t=j-d_{vH}}^{j-d\ v(j)-1}\eta^v(i,t)^T W_{2v}\eta^v(i,t)\\
&-\alpha^{d_{vH}}\sum_{t=j-d_v(j)}^{j-d_{vL}-1}\eta^v(i,t)^T W_{2v}\eta^v(i,t)
\end{aligned}$$

Then the following equations hold for any matrices $M_h=\begin{bmatrix}M_1^h\\ M_2^h\end{bmatrix}$, $N_h=\begin{bmatrix}N_1^h\\ N_2^h\end{bmatrix}$, $S_h=\begin{bmatrix}S_1^h\\ S_2^h\end{bmatrix}$, $M_v=\begin{bmatrix}M_1^v\\ M_2^v\end{bmatrix}$, $N_v=\begin{bmatrix}N_1^v\\ N_2^v\end{bmatrix}$ and $S_v=\begin{bmatrix}S_1^v\\ S_2^v\end{bmatrix}$ with appropriate dimensions:

$$\begin{aligned}
0=2\alpha^{d_{hL}}&\left[x^h(i,j)^T M_1^h+x^h\left(i-d_h(i),j\right)^T M_2^h\right]\\
&\times\left[x^h(i,j)-x^h\left(i-d_{hL},j\right)-\sum_{r=i-d_{hL}}^{i-1}\eta^h(r,j)\right]
\end{aligned}\tag{21}$$

$$\begin{aligned}
0=2\alpha^{d_{hH}}&\left[x^h(i,j)^T N_1^h+x^h\left(i-d_h(i),j\right)^T N_2^h\right]\\
&\times\left[x^h\left(i-d_h(i),j\right)-x^h\left(i-d_{hH},j\right)-\sum_{r=i-d_{hH}}^{i-1-d_h(i)}\eta^h(r,j)\right]
\end{aligned}\tag{22}$$

$$\begin{aligned}
0=2\alpha^{d_{hH}}&\left[x^h(i,j)^T S_1^h+x^h\left(i-d_h(i),j\right)^T S_2^h\right]\\
&\times\left[x^h\left(i-d_{hL},j\right)-x^h\left(i-d_h(i),j\right)-\sum_{r=i-d_h(i)}^{i-1-d_{hL}}\eta^h(r,j)\right]
\end{aligned}\tag{23}$$

$$\begin{aligned}
0=2\alpha^{d_{vL}}&\left[x^v(i,j)^T M_1^v+x^v\left(i,j-d_v(j)\right)^T M_2^v\right]\\
&\times\left[x^v(i,j)-x^v\left(i,j-d_{vL}\right)-\sum_{t=i-d_{vL}}^{j-1}\eta^v(i,t)\right]
\end{aligned}\tag{24}$$

$$\begin{aligned}
0=2\alpha^{d_{vH}}&\left[x^v(i,j)^T N_1^v+x^v\left(i,j-d_v(j)\right)^T N_2^v\right]\\
&\times\left[x^v\left(i,j-d_v(j)\right)-x^v\left(i,j-d_{vH}\right)-\sum_{t=i-d_{vH}}^{j-1-d_v(j)}\eta^v(i,t)\right]
\end{aligned}\tag{25}$$

$$0 = 2\alpha^{d_{vH}} \left[ x^v(i,j)^T S_1^v + x^h(i, j-d_v(j))^T S_2^v \right]$$

$$\times \left[ x^v(i, j-d_{vL}) - x^v(i, j-d_v(j)) - \sum_{t=j-d_v(j)}^{j-1-d_{vL}} \eta^v(i,t) \right] \tag{26}$$

On the other hand, for any matrices $X_h = \begin{bmatrix} X_{11h} & X_{12h} \\ * & X_{22h} \end{bmatrix} > 0$ , $Y_h = \begin{bmatrix} Y_{11h} & Y_{12h} \\ * & Y_{22h} \end{bmatrix} > 0$ , $X_v = \begin{bmatrix} X_{11v} & X_{12v} \\ * & X_{22v} \end{bmatrix} > 0$, $Y_v = \begin{bmatrix} Y_{11v} & Y_{12v} \\ * & Y_{22v} \end{bmatrix} > 0$, the following equations also hold

$$0 = \alpha^{d_{hL}} \left( \sum_{r=i-d_{hL}}^{i-1} \xi_1^h(i,j)^T X_h \xi_1^h(i,j) - \sum_{r=i-d_{hL}}^{i-1} \xi_1^h(i,j)^T X_h \xi_1^h(i,j) \right)$$

$$= \alpha^{d_{hL}} \left( d_{hL} \xi_1^h(i,j)^T X_h \xi_1^h(i,j) - \sum_{r=i-d_{hL}}^{i-1} \xi_1^h(i,j)^T X_h \xi_1^h(i,j) \right) \tag{27}$$

$$0 = \alpha^{d_{hH}} \left( \sum_{r=i-d_{hH}}^{i-1-d_{hL}} \xi_1^h(i,j)^T Y_h \xi_1^h(i,j) - \sum_{r=i-d_{hH}}^{i-1-d_{hL}} \xi_1^h(i,j)^T Y_h \xi_1^h(i,j) \right)$$

$$= \alpha^{d_{hH}} \left( (d_{hH} - d_{hL}) \xi_1^h(i,j)^T Y_h \xi_1^h(i,j) - \sum_{r=i-d_h(i)}^{i-1-d_{hL}} \xi_1^h(i,j)^T Y_h \xi_1^h(i,j) \right.$$

$$\left. - \sum_{r=i-d_{hH}}^{i-1-d_h(i)} \xi_1^h(i,j)^T Y_h \xi_1^h(i,j) \right) \tag{28}$$

$$0 = \alpha^{d_{vL}} \left( \sum_{t=j-d_{vL}}^{j-1} \xi_1^v(i,j)^T X_v \xi_1^v(i,j) - \sum_{t=j-d_{vL}}^{j-1} \xi_1^v(i,j)^T X_v \xi_1^v(i,j) \right)$$

$$= \alpha^{d_{vL}} \left( d_{vL} \xi_1^v(i,j)^T X_v \xi_1^v(i,j) - \sum_{t=j-d_{vL}}^{j-1} \xi_1^v(i,j)^T X_v \xi_1^v(i,j) \right) \tag{29}$$

$$0 = \alpha^{d_{vH}} \left( \sum_{t=j-d_{vH}}^{j-1-d_{vL}} \xi_1^v(i,j)^T Y_v \xi_1^v(i,j) - \sum_{t=j-d_{vH}}^{j-1-d_{vL}} \xi_1^v(i,j)^T Y_v \xi_1^v(i,j) \right)$$

$$= \alpha^{d_{vH}} \left( (d_{vH} - d_{vL}) \xi_1^v(i,j)^T Y_v \xi_1^v(i,j) - \sum_{t=j-d_v(i)}^{j-1-d_{vL}} \xi_1^v(i,j)^T Y_v \xi_1^v(i,j) \right.$$

$$\left. - \sum_{t=j-d_{vH}}^{j-1-d_v(j)} \xi_1^v(i,j)^T Y_v \xi_1^v(i,j) \right) \tag{30}$$

where

$$\xi_1^h(i,j) = \left[ x^h(i,j) \;\; x^h(i-d_h(i), j) \right]^T, \;\; \xi_1^v(i,j) = \left[ x^v(i,j) \;\; x^v(i, j-d_v(j)) \right]^T.$$

Combining (20)-(30) yields

$$V^h\left(x^h(i+1,j)\right)-\alpha V^h\left(x^h(i,j)\right)+V^v\left(x^v(i,j+1)\right)-\alpha V^v\left(x^v(i,j)\right)$$
$$=\xi(i,j)^T\Phi\xi(i,j)-\sum_{r=i-d_{hL}}^{i-1}\xi_2^h(r,j)^T\Psi_{1h}\xi_2^h(r,j)$$
$$-\sum_{r=i-d_h(i)}^{i-d_{hL}-1}\xi_2^h(r,j)^T\Psi_{2h}\xi_2^h(r,j)-\sum_{r=i-d_{hH}}^{i-d_h(i)-1}\xi_2^h(r,j)^T\Psi_{3h}\xi_2^h(r,j)$$
$$-\sum_{t=j-d_{vL}}^{j-1}\xi_2^v(i,t)^T\Psi_{1v}\xi_2^v(i,t)-\sum_{t=j-d_v(j)}^{j-d_{vL}-1}\xi_2^v(i,t)^T\Psi_{2v}\xi_2^v(i,t)$$
$$-\sum_{t=j-d_{vH}}^{j-d_v(j)-1}\xi_2^v(i,t)^T\Psi_{3v}\xi_2^v(i,t)$$

where

$$\xi_2^h(r,j)=\begin{bmatrix}\xi_1^h(r,j)^T & \eta^h(r,j)^T\end{bmatrix}^T,\ \ \xi_2^v(i,t)=\begin{bmatrix}\xi_1^v(i,t)^T & \eta^v(i,t)^T\end{bmatrix}^T,$$
$$\xi(i,j)=\begin{bmatrix}x(i,j) & x_d(i,j) & x_{dL}(i,j) & x_{dH}(i,j)\end{bmatrix}^T.$$

with

$$x(i,j)=\begin{bmatrix}x^h(i,j)\\ x^v(i,j)\end{bmatrix},\ \ x_d(i,j)=\begin{bmatrix}x^h(i-d(i),j)\\ x^v(i,j-d(j))\end{bmatrix},$$
$$x_{dL}(i,j)=\begin{bmatrix}x^h(i-d_{hL},j)\\ x^v(i,j-d_{vL})\end{bmatrix},\ \ x_{dH}(i,j)=\begin{bmatrix}x^h(i-d_{hH},j)\\ x^v(i,j-d_{vH})\end{bmatrix}.$$

Therefore, if inequalities (16)-(17) hold, then we can get that

$$V^h\left(x^h(i+1,j)\right)-\alpha V^h\left(x^h(i,j)\right)+V^v\left(x^v(i,j+1)\right)-\alpha V^v\left(x^v(i,j)\right)<0 \tag{31}$$

For simplicity, we define

$$V^h(i,j)=V^h\left(x^h(i,j)\right),\ \ V^v(i,j)=V^v\left(x^v(i,j)\right),\ \ V(i,j)=V\left(x(i,j)\right)$$

Thus, it is easy to get that

$$V^h(i+1,j)+V^v(i,j+1)<\alpha\left(V^h(i,j)+V^v(i,j)\right) \tag{32}$$

Since for any nonnegative integer $D>z=\max(z_1,z_2)$, we have $V^h(0,D)=V^v(D,0)=0$, then summing up both sides of (32) from $D-1$ to 0 with respect to $j$ and 0 to $D-1$ with respect to $i$, one gets

$$\sum_{i+j=D}V(i,j)=\ V^h(0,D)+V^h(1,D-1)+V^h(2,D-2)+\cdots+V^h(D-1,1)+V^h(D,0)$$

$$+V^v(0,D)+V^v(1,D-1)+V^v(2,D-2)+\cdots+V^v(D-1,1)+V^v(D,0)$$

$$<\alpha\Big(V^h(0,D-1)+V^v(0,D-1)+V^h(1,D-2)+V^v(1,D-2)$$

$$+\cdots+V^h(D-1,0)+V^v(D-1,0)\Big)$$

$$=\alpha\sum_{i+j=D-1}V(i,j)<\ldots<\alpha^{D-z}\sum_{i+j=z}V(i,j) \tag{33}$$

The proof is completed.

**Remark 3** Lemmas 4 provides the method for the estimation of Lyapunov functional candidate which will be used to design the controller for the 2D discrete switched system.

## 3.2. Robust reliable controller design

Consider systems (1), under the reliable controller $u^f(i,j)=\Omega^{\sigma(m)}K^{\sigma(m)}x(i,j)$, the corresponding closed-loop system with the bounded condition (11) is given by

$$\begin{bmatrix} x^h(i+1,j) \\ x^v(i,j+1) \end{bmatrix}=\left(\hat{A}^{\sigma(m)}+B^{\sigma(m)}\Omega^{\sigma(m)}K^{\sigma(m)}\right)\begin{bmatrix} x^h(i,j) \\ x^v(i,j) \end{bmatrix}+A_d^{\sigma(m)}\begin{bmatrix} x^h(i-d_h(i),j) \\ x^v(i,j-d_v(j)) \end{bmatrix} \tag{34}$$

Suppose that the $k$ th subsystem is activated at the switching instant $m_\kappa$, the $l$ th subsystem is activated at the switching instant $m_{\kappa+1}$.

**Theorem 1** Consider system (34), for given positive scalars $\alpha<1$, $\delta_k$, $\varepsilon_k$, if there exist positive definite symmetric matrices $\bar{P}^k=diag\{\bar{P}_h^k,\bar{P}_v^k\}$, $\bar{Q}_1^k=diag\{\bar{Q}_{1h}^k,\bar{Q}_{1v}^k\}$, $\bar{Q}_2^k=diag\{\bar{Q}_{2h}^k,\bar{Q}_{2v}^k\}$, $\bar{W}_1^k=diag\{\bar{W}_{1h}^k,\bar{W}_{1v}^k\}$, $W_2^k=diag\{\bar{W}_{2h}^k,\bar{W}_{2v}^k\}$, $\bar{X}_h^k=\begin{bmatrix}\bar{X}_{11h}^k & \bar{X}_{12h}^k \\ * & \bar{X}_{22h}^k\end{bmatrix}$, $\bar{X}_v^k=\begin{bmatrix}\bar{X}_{11v}^k & \bar{X}_{12v}^k \\ * & \bar{X}_{22v}^k\end{bmatrix}$, $\bar{Y}_h^k=\begin{bmatrix}\bar{Y}_{11h}^k & \bar{Y}_{12h}^k \\ * & \bar{Y}_{22h}^k\end{bmatrix}$, $\bar{Y}_v^k=\begin{bmatrix}\bar{Y}_{11v}^k & \bar{Y}_{12v}^k \\ * & \bar{Y}_{22v}^k\end{bmatrix}$, and any matrices $\bar{M}_h^k=\begin{bmatrix}\bar{M}_{1h}^k \\ \bar{M}_{2h}^k\end{bmatrix}$, $\bar{N}_h^k=\begin{bmatrix}\bar{N}_{1h}^k \\ \bar{N}_{2h}^k\end{bmatrix}$, $\bar{S}_h^k=\begin{bmatrix}\bar{S}_{1h}^k \\ \bar{S}_{2h}^k\end{bmatrix}$, $\bar{M}_v^k=\begin{bmatrix}\bar{M}_{1v}^k \\ \bar{M}_{2v}^k\end{bmatrix}$, $\bar{N}_v^k=\begin{bmatrix}\bar{N}_{1v}^k \\ \bar{N}_{2v}^k\end{bmatrix}$, $\bar{S}_v^k=\begin{bmatrix}\bar{S}_{1v}^k \\ \bar{S}_{2v}^k\end{bmatrix}$, $\Upsilon^k$ and $Z=diag\{Z_h,Z_v\}$ with appropriate dimensions, $k\in\underline{N}$, such that

$$\bar{\Phi} = \begin{bmatrix} \bar{\Phi}_1 & \bar{\Phi}_2^T & \bar{\Phi}_2^T & \bar{\Phi}_3^T & \mathrm{T}^k & \bar{\Phi}_7^T \\ * & \bar{\Phi}_4 & \bar{\Phi}_6 & \bar{\Phi}_6 & 0 & 0 \\ * & * & \bar{\Phi}_5 & \bar{\Phi}_6 & 0 & 0 \\ * & * & * & \Pi^k & 0 & 0 \\ * & * & * & * & -\delta_k I & 0 \\ * & * & * & * & * & -\varepsilon_k \Xi^k \end{bmatrix} < 0 \tag{35}$$

$$\bar{\Psi}^k = diag\left\{\bar{\Psi}_{1h}^k, \bar{\Psi}_{2h}^k, \bar{\Psi}_{3h}^k, \bar{\Psi}_{1v}^k, \bar{\Psi}_{2v}^k, \bar{\Psi}_{3v}^k\right\} \geq 0 \tag{36}$$

where

$\Lambda_1 = diag\left\{d_{hL} I_h, d_{vL} I_v\right\}$, $\Lambda_2 = diag\left\{\left(d_{hH} - d_{hL}\right) I_h, \left(d_{vH} - d_{vL}\right) I_v\right\}$,

$\Lambda_3 = diag\left\{\alpha^{d_{hL}} I_h, \alpha^{d_{vL}} I_v\right\}$, $\Lambda_4 = diag\left\{\alpha^{d_{hH}} I_h, \alpha^{d_{vH}} I_v\right\}$.

$$\bar{\Phi}_1 = \begin{bmatrix} \bar{\Phi}_{11} & \bar{\Phi}_{12} & \bar{\Phi}_{13} & \bar{\Phi}_{14} \\ * & \bar{\Phi}_{22} & \bar{\Phi}_{23} & \bar{\Phi}_{24} \\ * & * & \bar{\Phi}_{33} & 0 \\ * & * & * & \bar{\Phi}_{44} \end{bmatrix}$$

$\bar{\Phi}_{11} = \bar{Q}_1^k + \bar{Q}_2^k + \Lambda_3\left(\bar{M}_1^k + \bar{M}_1^{kT} + \Lambda_1 \bar{X}_{11}^k\right) + \Lambda_2 \Lambda_4 \bar{Y}_{11}^k - \alpha \bar{P}^k$, $\bar{\Phi}_{13} = -\Lambda_3 \bar{M}_1^k + \Lambda_4 \bar{S}_1^k$,

$\bar{\Phi}_{12} = \Lambda_3\left(\Lambda_1 \bar{X}_{12}^k + \bar{M}_{2k}^{kT}\right) + \Lambda_4\left(\bar{N}_1^k - \bar{S}_1^k + \Lambda_2 \bar{Y}_{12}^k\right)$, $\bar{\Phi}_{23} = -\Lambda_3 \bar{M}_2^k + \Lambda_4 \bar{S}_2^k$,

$\bar{\Phi}_{14} = -\Lambda_4 \bar{N}_1^k$, $\bar{\Phi}_{22} = \Lambda_1 \Lambda_3 \bar{X}_{22}^k + \Lambda_4\left(\bar{N}_2^k + \bar{N}_2^{kT} - \bar{S}_2^k - \bar{S}_2^{kT} + \Lambda_2 \bar{Y}_{22}^k\right)$, $\bar{\Phi}_{24} = -\Lambda_4 N_2^k$,

$\bar{\Phi}_{33} = -\Lambda_3 \bar{Q}_1^k$, $\bar{\Phi}_{44} = -\Lambda_4 \bar{Q}_2^k$, $\bar{\Phi}_2 = \left[\left(A^k - I\right) Z + B^k \Omega_0^k \Upsilon^k \quad A_d^k Z \quad 0 \quad 0\right]$,

$\bar{\Phi}_4 = -\left(\Lambda_1\right)^{-1} \bar{W}_1^k + \delta_k H^k H^{kT} + \varepsilon_k B^k \Omega_0^k \Xi^k \Omega_0^{kT} B^{kT}$,

$\bar{\Phi}_3 = \left[A^k Z + B^k \Omega_0^k \Upsilon^k \quad A_d^k Z \quad 0 \quad 0\right]$, $\bar{\Phi}_6 = \delta_k H^k H^{kT} + \varepsilon_k B^k \Omega_0^k \Xi^k \Omega_0^{kT} B^{kT}$,

$\bar{\Phi}_5 = -\left(\Lambda_2\right)^{-1} \bar{W}_2^k + \delta_k H^k H^{kT} + \varepsilon_k B^k \Omega_0^k \Xi^k \Omega_0^{kT} B^{kT}$, $\mathrm{T}^k = \left[\left(E^k Z\right)^T \quad \left(E_d^k Z\right)^T \quad 0 \quad 0\right]$,

$\Pi^k = -Z - Z^T + \bar{P}^k + \delta_k H^k H^{kT} + \varepsilon_k B^k \Omega_0^k \Xi^k \Omega_0^{kT} B^{kT}$, $\bar{\Phi}_7 = \left[\Xi^k \Upsilon^k \quad 0 \quad 0 \quad 0\right]$.

$\bar{X}_{11}^k = diag\left\{\bar{X}_{11h}^k, \bar{X}_{11v}^k\right\}$, $\bar{X}_{12}^k = diag\left\{\bar{X}_{12h}^k, \bar{X}_{12v}^k\right\}$, $\bar{X}_{22}^k = diag\left\{\bar{X}_{22h}^k, \bar{X}_{22v}^k\right\}$,

$\bar{Y}_{11}^k = diag\left\{\bar{Y}_{11h}^k, \bar{Y}_{11v}^k\right\}$, $\bar{Y}_{12}^k = diag\left\{\bar{Y}_{12h}^k, \bar{Y}_{12v}^k\right\}$, $\bar{Y}_{22}^k = diag\left\{\bar{Y}_{22h}^k, \bar{Y}_{22v}^k\right\}$,

$\bar{M}_1^k = diag\left\{\bar{M}_{1h}^k, \bar{M}_{1v}^k\right\}$, $\bar{M}_2^k = diag\left\{\bar{M}_{2h}^k, \bar{M}_{2v}^k\right\}$, $\bar{N}_1^k = diag\left\{\bar{N}_{1h}^k, \bar{N}_{1v}^k\right\}$,

$\bar{N}_2^k = diag\left\{\bar{N}_{2h}^k, \bar{N}_{2v}^k\right\}$, $\bar{S}_1^k = diag\left\{\bar{S}_{1h}^k, \bar{S}_{1v}^k\right\}$, $\bar{S}_2^k = diag\left\{\bar{S}_{2h}^k, \bar{S}_{2v}^k\right\}$.

$$\bar{\Psi}_{1h}^{k}=\begin{bmatrix}\bar{X}_h^k & \bar{M}_h^k \\ * & Z_h+Z_h^T-\bar{W}_{1h}^k\end{bmatrix},\ \bar{\Psi}_{2h}^{k}=\begin{bmatrix}\bar{Y}_h^k & \bar{N}_h^k \\ * & Z_h+Z_h^T-\bar{W}_{2h}^k\end{bmatrix},\ \bar{\Psi}_{3h}^{k}=\begin{bmatrix}\bar{Y}_h^k & \bar{S}_h^k \\ * & Z_h+Z_h^T-\bar{W}_{2h}^k\end{bmatrix},$$

$$\bar{\Psi}_{1v}^{k}=\begin{bmatrix}\bar{X}_v^k & M_v^k \\ * & Z_v+Z_v^T-\bar{W}_{1v}^k\end{bmatrix},\ \bar{\Psi}_{2v}^{k}=\begin{bmatrix}\bar{Y}_v^k & \bar{N}_v^k \\ * & Z_v+Z_v^T-\bar{W}_{2v}^k\end{bmatrix},\ \bar{\Psi}_{3v}^{k}=\begin{bmatrix}\bar{Y}_v^k & \bar{S}_v^k \\ * & Z_v+Z_v^T-\bar{W}_{2v}^k\end{bmatrix}.$$

Then, under the reliable controller

$$u^f(i,j)=\Omega^{\sigma(m)}K^{\sigma(m)}x(i,j),\ K^k=\Upsilon^k Z^{-1} \tag{37}$$

and the following average dwell time scheme

$$\tau_a>\tau_a^*=\frac{\ln\mu}{-\ln\alpha} \tag{38}$$

the corresponding closed-loop system is exponentially stable, where $\mu\geq 1$ satisfying

$$\bar{P}^k\leq\mu\bar{P}^l,\ \bar{P}^l\leq\mu\bar{P}^k,\ \bar{Q}_1^k\leq\mu\bar{Q}_1^l,\ \bar{Q}_1^l\leq\mu\bar{Q}_1^k,\ \bar{Q}_2^k\leq\mu\bar{Q}_2^l$$

$$\bar{Q}_2^l\leq\mu\bar{Q}_2^k,\ \bar{W}_1^k\leq\mu\bar{W}_1^l,\ \bar{W}_1^l\leq\mu\bar{W}_1^k,\ \bar{W}_2^k\leq\mu\bar{W}_2^l,\ \bar{W}_2^l\leq\mu\bar{W}_2^k \tag{39}$$

**Proof** We consider the following Lyapunov function candidate for the $k$ th$\left(k\in\underline{N}\right)$ subsystem

$$V_k\left(x(i,j)\right)=V_k^h\left(x^h(i,j)\right)+V_k^v\left(x^v(i,j)\right) \tag{40}$$

where

$$V_k^h\left(x^h(i,j)\right)=\sum_{g=1}^{3}V_{gk}^h\left(x^h(i,j)\right)$$

$$V_{1k}^h\left(x^h(i,j)\right)=x^h(i,j)^T P_h^k x^h(i,j)$$

$$V_{2k}^h\left(x^h(i,j)\right)=\sum_{r=i-d_{hL}}^{i-1}x^h(r,j)^T Q_{1h}^k x^h(r,j)\alpha^{i-r-1}+\sum_{r=i-d_{hH}}^{i-1}x^h(r,j)^T Q_{2h}^k x^h(r,j)\alpha^{i-r-1}$$

$$V_{3k}^h\left(x^h(i,j)\right)=\sum_{s=-d_{hL}}^{-1}\sum_{r=i+s}^{i-1}\eta^h(r,j)^T W_{1h}^k\eta^h(r,j)\alpha^{i-r-1}+\sum_{s=-d_{hH}}^{-d_{hL}-1}\sum_{r=i+s}^{i-1}\eta^h(r,j)^T W_{2h}^k\eta^h(r,j)\alpha^{i-r-1}$$

$$V_k^v\left(x^v(i,j)\right)=\sum_{g=1}^{3}V_{gk}^v\left(x^v(i,j)\right)$$

$$V_{1k}^v\left(x^v(i,j)\right)=x^v(i,j)^T P_v^k x^v(i,j)$$

$$V_{2k}^{v}\left(x^{v}\left(i,j\right)\right)=\sum_{t=j-d_{vL}}^{j-1} x^{v}\left(i,t\right)^{T} Q_{1v}^{k} x^{v}\left(i,t\right)\alpha^{j-t-1}$$
$$+\sum_{t=j-d_{vH}}^{j-1} x^{v}\left(i,t\right)^{T} Q_{2v}^{k} x^{v}\left(i,t\right)\alpha^{j-t-1}$$

$$V_{3k}^{v}\left(x^{v}\left(i,j\right)\right)=\sum_{s=-d_{vL}}^{-1}\sum_{t=j+s}^{j-1}\eta^{v}\left(i,t\right)^{T} W_{1v}^{k}\eta^{v}\left(i,t\right)\alpha^{j-t-1}$$
$$+\sum_{s=-d_{vH}}^{-d_{vL}-1}\sum_{t=j+s}^{j-1}\eta^{v}\left(i,t\right)^{T} W_{2v}^{k}\eta^{v}\left(i,t\right)\alpha^{j-t-1}$$

By Lemma 4, we know that for a given scalar $\alpha<1$, if there exist positive definite symmetric matrices $P^{k}=diag\left\{P_{h}^{k},P_{v}^{k}\right\}$, $Q_{1}^{k}=diag\left\{Q_{1h}^{k},Q_{1v}^{k}\right\}$, $Q_{2}^{k}=diag\left\{Q_{2h}^{k},Q_{2v}^{k}\right\}$, $W_{1}^{k}=diag\left\{W_{1h}^{k},W_{1v}^{k}\right\}$, $W_{2}^{k}=diag\left\{W_{2h}^{k},W_{2v}^{k}\right\}$, $X_{h}^{k}=\begin{bmatrix} X_{11h}^{k} & X_{12h}^{k} \\ * & X_{22h}^{k}\end{bmatrix}$, $X_{v}^{k}=\begin{bmatrix} X_{11v}^{k} & X_{12v}^{k} \\ * & X_{22v}^{k}\end{bmatrix}$, $Y_{h}^{k}=\begin{bmatrix} Y_{11h}^{k} & Y_{12h}^{k} \\ * & Y_{22h}^{k}\end{bmatrix}$, $Y_{v}^{k}=\begin{bmatrix} Y_{11v}^{k} & Y_{12v}^{k} \\ * & Y_{22v}^{k}\end{bmatrix}$, $M_{h}^{k}=\begin{bmatrix} M_{1h}^{k} \\ M_{2h}^{k}\end{bmatrix}$, $N_{h}^{k}=\begin{bmatrix} N_{1h}^{k} \\ N_{2h}^{k}\end{bmatrix}$, $S_{h}^{k}=\begin{bmatrix} S_{1h}^{k} \\ S_{2h}^{k}\end{bmatrix}$, $M_{v}^{k}=\begin{bmatrix} M_{1v}^{k} \\ M_{2v}^{k}\end{bmatrix}$, $N_{v}^{k}=\begin{bmatrix} N_{1v}^{k} \\ N_{2v}^{k}\end{bmatrix}$, and $S_{v}^{k}=\begin{bmatrix} S_{1v}^{k} \\ S_{2v}^{k}\end{bmatrix}$ with appropriate dimensions, such that

$$\hat{\Phi}=\begin{bmatrix} \hat{\Phi}_{1} & \hat{\Phi}_{2}^{T} & \hat{\Phi}_{2}^{T} & \hat{\Phi}_{3}^{T}P^{k} \\ * & -\left(\Lambda_{1}\right)^{-1}W_{1k}^{-1} & 0 & 0 \\ * & * & -\left(\Lambda_{2}\right)^{-2}W_{2k}^{-1} & 0 \\ * & * & * & -P^{k}\end{bmatrix}<0 \tag{41}$$

$$\Psi^{k}=diag\left\{\Psi_{1h}^{k},\Psi_{2h}^{k},\Psi_{3h}^{k},\Psi_{1v}^{k},\Psi_{2v}^{k},\Psi_{3v}^{k}\right\}\geq 0 \tag{42}$$

where

$$\hat{\Phi}_{1}=\begin{bmatrix} \hat{\Phi}_{11} & \hat{\Phi}_{12} & \hat{\Phi}_{13} & \hat{\Phi}_{14} \\ * & \hat{\Phi}_{22} & \hat{\Phi}_{23} & \hat{\Phi}_{24} \\ * & * & \hat{\Phi}_{33} & 0 \\ * & * & * & \hat{\Phi}_{44}\end{bmatrix}$$

$\hat{\Phi}_{11}=Q_{1}^{k}+Q_{2}^{k}+\Lambda_{3}\left(M_{1}^{k}+M_{1}^{kT}+\Lambda_{1}X_{11}^{k}\right)+\Lambda_{2}\Lambda_{4}Y_{11}^{k}-\alpha P^{k}$, $\hat{\Phi}_{13}=-\Lambda_{3}M_{1}^{k}+\Lambda_{4}S_{1}^{k}$,

$\hat{\Phi}_{12}=\Lambda_{3}\left(\Lambda_{1}X_{12}^{k}+M_{2}^{kT}\right)+\Lambda_{4}\left(N_{1}^{k}-S_{1}^{k}+\Lambda_{2}Y_{12}^{k}\right)$, $\hat{\Phi}_{14}=-\Lambda_{4}N_{1}^{k}$,

$\hat{\Phi}_{22}=\Lambda_{1}\Lambda_{3}X_{22}^{k}+\Lambda_{4}\left(N_{2}^{k}+N_{2}^{kT}-S_{2}^{k}-S_{2}^{kT}+\Lambda_{2}Y_{22}^{k}\right)$, $\hat{\Phi}_{23}=-\Lambda_{3}M_{2}^{k}+\Lambda_{4}S_{2}^{k}$,

$\hat{\Phi}_{24} = -\Lambda_4 N_2^k$, $\hat{\Phi}_{33} = -\Lambda_3 Q_1^k$, $\hat{\Phi}_{44} = -\Lambda_4 Q_2^k$,

$\hat{\Phi}_2 = \left[ \hat{A}^k + B^k \Omega^k K^k - I \quad \hat{A}_d^k \quad 0 \quad 0 \right]$, $\hat{\Phi}_3 = \left[ \hat{A}^k + B^k \Omega^k K^k \quad \hat{A}_d^k \quad 0 \quad 0 \right]$.

$X_{11}^k = diag\left\{ X_{11h}^k, X_{11v}^k \right\}$, $X_{12}^k = diag\left\{ X_{12h}^k, X_{12v}^k \right\}$, $X_{22}^k = diag\left\{ X_{22h}^k, X_{22v}^k \right\}$

$Y_{11}^k = diag\left\{ Y_{11h}^k, Y_{11v}^k \right\}$, $Y_{12}^k = diag\left\{ Y_{12h}^k, Y_{12v}^k \right\}$, $Y_{22}^k = diag\left\{ Y_{22h}^k, Y_{22v}^k \right\}$,

$M_1^k = diag\left\{ M_{1h}^k, M_{1v}^k \right\}$, $M_2^k = diag\left\{ M_{2h}^k, M_{2v}^k \right\}$, $N_1^k = diag\left\{ N_{1h}^k, N_{1v}^k \right\}$,

$N_2^k = diag\left\{ N_{2h}^k, N_{2v}^k \right\}$, $S_1^k = diag\left\{ S_{1h}^k, S_{1v}^k \right\}$, $S_2^k = diag\left\{ S_{2h}^k, S_{2v}^k \right\}$.

$$\Psi_{1h}^k = \begin{bmatrix} X_h^k & M_h^k \\ * & W_{1h}^k \end{bmatrix},\ \Psi_{2h}^k = \begin{bmatrix} Y_h^k & N_h^k \\ * & W_{2h}^k \end{bmatrix},\ \Psi_{3h}^k = \begin{bmatrix} Y_h^k & S_h^k \\ * & W_{2h}^k \end{bmatrix},$$

$$\Psi_{1v}^k = \begin{bmatrix} X_v^k & M_v^k \\ * & W_{1v}^k \end{bmatrix},\ \Psi_{2v}^k = \begin{bmatrix} Y_v^k & N_v^k \\ * & W_{2v}^k \end{bmatrix},\ \Psi_{3v}^k = \begin{bmatrix} Y_v^k & S_v^k \\ * & W_{2v}^k \end{bmatrix}.$$

then the following inequality holds

$$\sum_{i+j=D} V_k(i,j) < \alpha^{D-D^k} \sum_{i+j=D^k} V_k(i,j) \tag{43}$$

where $D^k \in [z, D)$.

Now Let $\chi = N_{\sigma(m)}(z, D)$ denote the switch number of $\sigma(m)$ on an interval $(z, D)$, and let $m_{\kappa-\chi+1} < m_{\kappa-\chi+2} < \cdots < m_{\kappa-1} < m_\kappa$ denote the switching points of $\sigma(m)$ over the interval $(z, D)$. Using the condition (39), (40), then at switching instant $m_\kappa = i + j$, we can get the following inequality

$$\sum_{i+j=m_\kappa} V_{\sigma(m_\kappa)}(i,j) \le \mu \sum_{i+j=m_\kappa} V_{\sigma(m_{\kappa-1})}(i,j) \tag{44}$$

Thus, it follows that

$$\begin{aligned} &\sum_{i+j=D} V_{\sigma(m_\kappa)}(i,j) < \alpha^{D-m_k} \sum_{i+j=m_\kappa} V_{\sigma(m_\kappa)}(i,j) \\ &\le \mu\alpha^{D-m_\kappa} \sum_{i+j=m_\kappa} V_{\sigma(m_{\kappa-1})}(i,j)^- \\ &< \mu^2 \alpha^{D-m_\kappa} \alpha^{m_\kappa - m_{\kappa-1}} \sum_{i+j=m_{\kappa-1}} V_{\sigma(m_{\kappa-2})}(i,j)^- \\ &= \mu^2 \alpha^{D-m_{\kappa-1}} \sum_{i+j=m_{\kappa-1}} V_{\sigma(m_{\kappa-2})}(i,j)^- < \ldots \\ &< \mu^\chi \alpha^{D-z} \sum_{i+j=z} V_{\sigma(m_{\kappa-\chi+1})}(i,j) \end{aligned} \tag{45}$$

According to definition 2, we know

$$\chi = N_{\sigma(m)}\left(z, D\right) \le N_0 + \frac{D-z}{\tau_a} \tag{46}$$

Then when $0<\alpha<1$, from (45), we have

$$\begin{aligned}\sum_{i+j=D} V_{\sigma(m_\kappa)}\left(i,j\right) &< \mu^{\chi}\alpha^{D-z}\sum_{i+j=z} V_{\sigma\left(m_{k-\chi+1}\right)}\left(i,j\right)\\ &= \mu^{N_0} e^{-\left(-\frac{\ln\mu}{\tau_a}-\ln\alpha\right)(D-z)}\sum_{i+j=z} V_{\sigma\left(m_{\kappa-\chi+1}\right)}\left(i,j\right)\end{aligned} \tag{47}$$

Then the following inequality can be obtained

$$\sum_{i+j=D}\left\|x\left(i,j\right)\right\|^2 \le \left(\zeta_1/\zeta_2\right)\mu^{N_0} e^{-\left(-\frac{\ln\mu}{\tau_a}-\ln\alpha\right)(D-z)}\sum_{i+j=z}\left\|x\left(i,j\right)\right\|_C^2 \tag{48}$$

where

$$\zeta_1 = \max_{k,l\in\underline{N},k\neq l}\left\{\lambda_1,\lambda_2\right\}$$

$$\begin{aligned}\lambda_1 = {}&\lambda_{\max}\left(P^k\right)+\max\left(d_{hL},d_{vL}\right)\lambda_{\max}\left(Q_1^k\right)+\max\left(d_{hH},d_{vH}\right)\lambda_{\max}\left(Q_2^k\right)\\ &+\max\left(d_{hL}^2,d_{vL}^2\right)\lambda_{\max}\left(W_1^k\right)+\max\left(\left(d_{hH}-d_{hL}\right)^2,\left(d_{vH}-d_{vL}\right)^2\right)\lambda_{\max}\left(W_2^k\right)\end{aligned}$$

$$\begin{aligned}\lambda_2 = {}&\lambda_{\max}\left(P^l\right)+\max\left(d_{hL},d_{vL}\right)\lambda_{\max}\left(Q_1^l\right)+\max\left(d_{hH},d_{vH}\right)\lambda_{\max}\left(Q_2^l\right)\\ &+\max\left(d_{hL}^2,d_{vL}^2\right)\lambda_{\max}\left(W_1^l\right)+\max\left(\left(d_{hH}-d_{hL}\right)^2,\left(d_{vH}-d_{vL}\right)^2\right)\lambda_{\max}\left(W_2^l\right)\end{aligned}$$

$$\zeta_2 = \min_{k,l\in\underline{N},k\neq l}\left\{\lambda_{\min}\left(P^k\right),\lambda_{\min}\left(P^l\right)\right\}$$

$$\begin{aligned}\sum_{i+j=z}\left\|x\left(i,j\right)\right\|_C^2 \triangleq \sup_{\substack{-d_{hH}<\theta_h\le 0,\\ -d_{vH}<\theta_v\le 0}}\sum_{i+j=z}\Big\{&\left\|x^h\left(i-\theta_h,j\right)\right\|^2+\left\|x^v\left(i,j-\theta_v\right)\right\|^2,\\ &\left\|\eta^h\left(i-\theta_h,j\right)\right\|^2+\left\|\eta^v\left(i,j-\theta_v\right)\right\|^2\Big\}\end{aligned}$$

where

$$\eta^h\left(i-\theta_h,j\right) = x^h\left(i-\theta_h+1,j\right)-x^h\left(i-\theta_h,j\right)$$

$$\eta^v\left(i,j-\theta_v\right) = x^v\left(i,j-\theta_v+1\right)-x^v\left(i,j-\theta_v\right)$$

So when the condition (38) is satisfied, $\sum_{i+j=D}\left\|x\left(i,j\right)\right\|^2$ will tend to zero exponentially, which means the corresponding closed-loop 2D discrete switched system is exponentially stable.

Since for $P^k$, there exits a reversible matrix $G = diag\left\{G_h,G_v\right\}$ such that $G^T P^{-k} G \geq G^T + G - P^k$ holds, using $diag\{I \quad I \quad I \quad G^T P^{-k}\}$ and $diag\{I \quad I \quad I \quad P^{-k} G\}$ to pre- and post-multiply the left of (41), respectively, then can obtain that if the following inequality satisfied, then (41) holds.

$$\hat{\Phi} = \begin{bmatrix} \hat{\Phi}_1 & \hat{\Phi}_2^T & \hat{\Phi}_2^T & \hat{\Phi}_3^T G \\ * & -(\Lambda_1)^{-1} W_{1k}^{-1} & 0 & 0 \\ * & * & -(\Lambda_2)^{-2} W_{2k}^{-1} & 0 \\ * & * & * & -G - G^T + P^k \end{bmatrix} < 0 \tag{49}$$

Denoting $Z = G^{-1}$, $\bar{W}_1^k = W_{1k}^{-1}$, $\bar{W}_2^k = W_{2k}^{-1}$ and using $diag\{Z^T, Z^T, Z^T, Z^T, I, I, Z^T\}$ and $diag\{Z, Z, Z, Z, I, I, Z\}$ to pre- and post-multiply the left of (49), respectively, and assigning the following new matrices

$\bar{P}^k = Z^T P^k Z$ $\bar{Q}_1^k = Z^T Q_1^k Z$, $\bar{Q}_2^k = Z^T Q_2^k Z$, $\bar{X}_{11}^k = Z^T X_{11}^k Z$,

$\bar{X}_{12}^k = Z^T X_{12}^k Z$, $\bar{X}_{22}^k = Z^T X_{22}^k Z$, $\bar{Y}_{11}^k = Z^T Y_{11}^k Z$, $\bar{Y}_{12}^k = Z^T Y_{12}^k Z$,

$\bar{Y}_{22}^k = Z^T Y_{22}^k Z$, $\bar{M}_1^k = Z^T M_1^k Z$, $\bar{M}_2^k = Z^T M_2^k Z$, $\bar{N}_1^k = Z^T N_1^k Z$,

$\bar{N}_2^k = Z^T N_2^k Z$, $\bar{S}_1^k = Z^T S_1^k Z$, $\bar{S}_2^k = Z^T S_2^k Z$, $\Upsilon^k = K^k Z$.

Furthermore, replacing $\hat{A}^k$, $\hat{A}_d^k$ and $\Omega^k$ with $A^k + H^k F^k(i,j) E^k$, $A_d^k + H^k F^k(i,j) E_d^k$ and $\Omega_0^k (I + \Theta^k)$, then applying lemma 1, 2 and 3, inequality (35) is directly obtained.

In addition, denoting $Z_h = G_h^{-1}$, $Z_v = G_v^{-1}$, $\bar{W}_{1h}^k = (W_{1h}^k)^{-1}$, $\bar{W}_{1v}^k = (W_{1v}^k)^{-1}$, $\bar{W}_{2h}^k = (W_{2h}^k)^{-1}$, $\bar{W}_{2v}^k = (W_{2v}^k)^{-1}$, then using $diag\{Z_h^T, Z_h^T, Z_h^T, Z_h^T, Z_h^T, Z_h^T, Z_v^T, Z_v^T, Z_v^T, Z_v^T, Z_v^T, Z_v^T\}$ and $diag\{Z_h, Z_h, Z_h, Z_h, Z_h, Z_h, Z_v, Z_v, Z_v, Z_v, Z_v, Z_v\}$ to pre- and post-multiply the left of (42), we can obtain

$$\hat{\Psi}^k = diag\{\hat{\Psi}_{1h}^k, \hat{\Psi}_{2h}^k, \hat{\Psi}_{3h}^k, \hat{\Psi}_{1v}^k, \hat{\Psi}_{2v}^k, \hat{\Psi}_{3v}^k\} \geq 0 \tag{50}$$

where

$$\hat{\Psi}_{1h}^k = \begin{bmatrix} \bar{X}_h^k & \bar{M}_h^k \\ * & Z_h^T W_{1h}^k Z \end{bmatrix}, \quad \hat{\Psi}_{2h}^k = \begin{bmatrix} \bar{Y}_h^k & \bar{N}_h^k \\ * & Z_h^T W_{2h}^k Z_h \end{bmatrix}, \quad \hat{\Psi}_{3h}^k = \begin{bmatrix} \bar{Y}_h^k & \bar{S}_h^k \\ * & Z_h^T W_{2h}^k Z_h \end{bmatrix},$$

$$\hat{\Psi}_{1v}^{k}=\begin{bmatrix}\bar{X}_{v}^{k} & \bar{M}_{v}^{k}\\ * & Z_{v}^{T}W_{1v}^{k}Z_{v}\end{bmatrix},\ \hat{\Psi}_{2v}^{k}=\begin{bmatrix}\bar{Y}_{v}^{k} & \bar{N}_{v}^{k}\\ * & Z_{v}^{T}W_{2v}^{k}Z_{v}\end{bmatrix},\ \hat{\Psi}_{3v}^{k}=\begin{bmatrix}\bar{Y}_{v}^{k} & \bar{S}_{v}^{k}\\ * & Z_{v}^{T}W_{2v}^{k}Z_{v}\end{bmatrix}.$$

$\bar{X}_h^k = Z_h^T X_h^k Z_h$, $\bar{X}_v^k = Z_v^T X_v^k Z_v$, $\bar{Y}_h^k = Z_h^T Y_h^k Z_h$, $\bar{Y}_v^k = Z_v^T Y_v^k Z_v$, $\bar{M}_h^k = Z_h^T M_h^k Z_h$,

$\bar{M}_v^k = Z_v^T M_v^k Z_v$, $\bar{N}_h^k = Z_h^T N_h^k Z_h$, $\bar{N}_v^k = Z_v^T N_v^k Z_v$, $\bar{S}_h^k = Z_h^T S_h^k Z_h$, $\bar{S}_v^k = Z_v^T S_v^k Z_v$.

then applying the following relations

$$Z_h^T W_{1h}^k Z_h \geq Z_h^T + Z_h - \bar{W}_{1h}^k,\ \ Z_h^T W_{2h}^k Z_h \geq Z_h^T + Z_h - \bar{W}_{2h}^k,$$

$$Z_v^T W_{1v}^k Z_v \geq Z_v^T + Z_v - \bar{W}_{1v}^k,\ \ Z_v^T W_{2v}^k Z_v \geq Z_v^T + Z_v - \bar{W}_{2v}^k$$

we can know that if (36) are satisfied, then (42) hold. The proof is completed.

**Remark 4** It is noticed that (35)-(36) are LMIs, we can firstly solve the LMIs to obtain the solutions of matrices $\bar{P}^k$, $\bar{Q}_1^k$, $\bar{Q}_2^k$, $\bar{W}_1^k$, $W_2^k$, $\bar{X}_h^k$, $\bar{X}_v^k$, $\bar{Y}_h^k$, $\bar{Y}_v^k$, $\bar{M}_h^k$, $\bar{N}_h^k$, $\bar{S}_h^k$, $\bar{M}_v^k$, $\bar{N}_v^k$, $\bar{S}_v^k$, $Z$ and $\Upsilon^k$. Then $K^k$ can be obtained by $K^k = \Upsilon^k Z^{-1}$.

**Remark 5** It is noticed that when $\mu = 1$ in $\tau_a > \tau_a^* = \dfrac{\ln \mu}{-\ln \alpha}$, (39) turns out to be $\bar{P}^k = \bar{P}^l$, $\bar{Q}_1^k = \bar{Q}_1^l$, $\bar{Q}_2^k = \bar{Q}_2^l$, $\bar{W}_1^k = \bar{W}_1^l$, $\bar{W}_2^k = \bar{W}_2^l$, $\forall k,l \in \underline{N}$. In this case, we have $\tau_a(m) > \tau_a^*(m) = 0$, which means that the switching signal can be arbitrary.

**Remark 6** Note that when $d_{hH} = d_{hL} = d_h$ and $d_{vH} = d_{vL} = d_v$, the system (34) is reduced to the following system

$$\begin{bmatrix}x^h(i+1,j)\\ x^v(i,j+1)\end{bmatrix}=\left(\hat{A}^{\sigma(m)}+B^{\sigma(m)}\Omega^{\sigma(m)}K^{\sigma(m)}\right)\begin{bmatrix}x^h(i,j)\\ x^v(i,j)\end{bmatrix}+A_d^{\sigma(m)}\begin{bmatrix}x^h(i-d_h,j)\\ x^v(i,j-d_v)\end{bmatrix} \quad (51)$$

**Theorem 2** Consider system (51), for given positive scalars $\delta_k$, $\varepsilon_k$, $\alpha < 1$, if there exist positive definite symmetric matrices $\bar{P}^k = diag\left\{\bar{P}_h^k, \bar{P}_v^k\right\}$, $\bar{Q}_1^k = diag\left\{\bar{Q}_{1h}^k, \bar{Q}_{1v}^k\right\}$, $\bar{W}_1^k = diag\left\{\bar{W}_{1h}^k, \bar{W}_{1v}^k\right\}$, $\bar{X}_h^k = \begin{bmatrix}\bar{X}_{11h}^k & \bar{X}_{12h}^k\\ * & \bar{X}_{22h}^k\end{bmatrix}$, $\bar{X}_v^k = \begin{bmatrix}\bar{X}_{11v}^k & \bar{X}_{12v}^k\\ * & \bar{X}_{22v}^k\end{bmatrix}$, and any matrices $\Upsilon^k$, $\bar{M}_h^k = \begin{bmatrix}\bar{M}_{1h}^k\\ \bar{M}_{2h}^k\end{bmatrix}$, $\bar{M}_v^k = \begin{bmatrix}\bar{M}_{1v}^k\\ \bar{M}_{2v}^k\end{bmatrix}$, $Z = \mathrm{diag}\left\{Z_h, Z_v\right\}$, with appropriate dimensions, $k \in \underline{N}$, such that

$$\tilde{\Phi}=\begin{bmatrix}\tilde{\Phi}_1 & \tilde{\Phi}_2 & \tilde{\Phi}_3 & \mathrm{T}^k & \Upsilon^{kT}\Xi^k\\ * & \tilde{\Phi}_4 & \tilde{\Phi}_5 & 0 & 0\\ * & * & \Pi & 0 & 0\\ * & * & * & -\delta_k I & 0\\ * & * & * & * & -\varepsilon_k\Xi^k\end{bmatrix}<0 \tag{52}$$

$$\bar{\Psi}_{1h}^k=\begin{bmatrix}\bar{X}_h^k & \bar{M}_h^k\\ * & Z_h+Z_h^T-\bar{W}_{1h}^k\end{bmatrix}\ge 0 \tag{53}$$

$$\bar{\Psi}_{1v}^k=\begin{bmatrix}\bar{X}_v^k & M_v^k\\ * & Z_v+Z_v^T-\bar{W}_{1v}^k\end{bmatrix}\ge 0 \tag{54}$$

where

$$\tilde{\Lambda}_1=diag\left\{d_hI_h,d_vI_v\right\},\ \ \tilde{\Lambda}_3=diag\left\{\alpha^{d_h}I_h,\alpha^{d_v}I_v\right\},\ \ \tilde{\Phi}_1=\begin{bmatrix}\tilde{\Phi}_{11} & \tilde{\Phi}_{12}\\ * & \tilde{\Phi}_{22}\end{bmatrix}.$$

$$\tilde{\Phi}_{11}=\bar{Q}_1^k+\tilde{\Lambda}_3\left(\bar{M}_1^k+\bar{M}_1^{kT}+\tilde{\Lambda}_1\bar{X}_{11}^k\right)-\alpha\bar{P}^k,\tilde{\Phi}_{12}=\tilde{\Lambda}_3\left(\tilde{\Lambda}_1\bar{X}_{12}^k+\bar{M}_{2k}^{kT}\right)-\tilde{\Lambda}_3\bar{M}_1^k,$$

$$\tilde{\Phi}_{22}=\tilde{\Lambda}_1\tilde{\Lambda}_3\bar{X}_{22}^k-\tilde{\Lambda}_3\bar{M}_2^k-\left(\tilde{\Lambda}_3\bar{M}_2^k\right)^T-\tilde{\Lambda}_3\bar{Q}_1^k,$$

$$\tilde{\Phi}_4=-\left(\tilde{\Lambda}_1\right)^{-1}\bar{W}_1^k+\delta_kH^kH^{kT}+\varepsilon_kB^k\Omega_0^k\Xi^k\Omega_0^{kT}B^{kT},$$

$$\tilde{\Phi}_2=\left[\left(A^k-I\right)Z+B^k\Omega_0^k\Upsilon^k \quad A_d^kZ\right],\ \ \tilde{\Phi}_3=\left[A^kZ+B^k\Omega_0^k\Upsilon^k \quad A_d^kZ\right],$$

$$\Pi^k=-Z-Z^T+\bar{P}^k+\delta_kH^kH^{kT}+\varepsilon_kB^k\Omega_0^k\Xi^k\Omega_0^{kT}B^{kT},$$

$$\mathrm{T}^k=\left[\left(E^kZ\right)^T \quad \left(E_d^kZ\right)^T\right],\tilde{\Phi}_5=\delta_kH^kH^{kT}+\varepsilon_kB^k\Omega_0^k\Xi^k\Omega_0^{kT}B^{kT},$$

$$\bar{X}_{11}^k=diag\left\{\bar{X}_{11h}^k,\bar{X}_{11v}^k\right\},\ \ \bar{X}_{12}^k=diag\left\{\bar{X}_{12h}^k,\bar{X}_{12v}^k\right\},\ \ \bar{X}_{22}^k=diag\left\{\bar{X}_{22h}^k,\bar{X}_{22v}^k\right\},$$

$$\bar{M}_1^k=diag\left\{\bar{M}_{1h}^k,\bar{M}_{1v}^k\right\},\ \ \bar{M}_2^k=diag\left\{\bar{M}_{2h}^k,\bar{M}_{2v}^k\right\}.$$

Then, under the following reliable controller

$$u^f\left(i,j\right)=\Omega^{\sigma(m)}K^{\sigma(m)}x\left(i,j\right),\ \ K^k=\Upsilon^kZ^{-1} \tag{55}$$

and the average dwell time scheme

$$\tau_a>\tau_a^*=\frac{\ln\mu}{-\ln\alpha} \tag{56}$$

the corresponding closed-loop system is exponentially stable, where $\mu\ge 1$ satisfying

$$\bar{P}^k\le\mu\bar{P}^l,\ \ \bar{P}^l\le\mu\bar{P}^k,\ \ \bar{Q}_1^k\le\mu\bar{Q}_1^l,$$

$$\bar{Q}_1^l \le \mu \bar{Q}_1^k,\ \bar{W}_1^k \le \mu \bar{W}_1^l,\ \bar{W}_1^l \le \mu \bar{W}_1^k \tag{57}$$

**Proof** We consider the following Lyapunov function candidate for the $k$ th$\left(k \in \underline{N}\right)$ subsystem

$$V_k\left(x(i,j)\right) = V_k^h\left(x^h(i,j)\right) + V_k^v\left(x^v(i,j)\right) \tag{58}$$

where

$$V_k^h\left(x^h(i,j)\right) = \sum_{g=1}^{3} V_{gk}^h\left(x^h(i,j)\right)$$

$$V_{1k}^h\left(x^h(i,j)\right) = x^h(i,j)^T P_h^k x^h(i,j)$$

$$V_{2k}^h\left(x^h(i,j)\right) = \sum_{r=i-d_h}^{i-1} x^h(r,j)^T Q_{1h}^k x^h(r,j)\alpha^{i-r-1}$$

$$V_{3k}^h\left(x^h(i,j)\right) = \sum_{s=-d_h}^{-1} \sum_{r=i+s}^{i-1} \eta^h(r,j)^T W_{1h}^k \eta^h(r,j)\alpha^{i-r-1}$$

$$V_k^v\left(x^v(i,j)\right) = \sum_{g=1}^{3} V_{gk}^v\left(x^v(i,j)\right)$$

$$V_{1k}^v\left(x^v(i,j)\right) = x^v(i,j)^T P_v^k x^v(i,j)$$

$$V_{2k}^v\left(x^v(i,j)\right) = \sum_{t=j-d_v}^{j-1} x^v(i,t)^T Q_{1v}^k x^v(i,t)\alpha^{j-t-1}$$

$$V_{3k}^v\left(x^v(i,j)\right) = \sum_{s=-d_v}^{-1} \sum_{t=j+s}^{j-1} \eta^v(i,t)^T W_{1v}^k \eta^v(i,t)\alpha^{j-t-1}$$

The remainder process can be followed by the same lines of the proof of Theorem 1, and we omit the details.

# 4. Numerical example

In this section, we present an example to illustrate the effectiveness of the proposed approach. Consider system (1) with parameters as follows

$$A_1 = \begin{bmatrix} 1 & 1.5 \\ 0.6 & 0.5 \end{bmatrix},\ A_{d1} = \begin{bmatrix} -0.2 & 0 \\ -0.1 & -0.1 \end{bmatrix},\ B_1 = \begin{bmatrix} -4.5 & 0 \\ 1 & -3 \end{bmatrix},\ H_1 = \begin{bmatrix} 0.1 & 0.15 \\ 0.1 & 0.2 \end{bmatrix}$$

$$E_1 = \begin{bmatrix} 0.1 & 0.15 \\ 0.1 & 0 \end{bmatrix},\ E_{d1} = \begin{bmatrix} 0.1 & 0 \\ 0.1 & 0.1 \end{bmatrix},\ F_1 = \begin{bmatrix} \sin\left(0.5\pi(i+j)\right) & 0 \\ 0 & \sin\left(0.5\pi(i+j)\right) \end{bmatrix}$$

$$A_2=\begin{bmatrix}1.5 & 1.5\\ 0.5 & 0.5\end{bmatrix},\ A_{d2}=\begin{bmatrix}-0.15 & 0.1\\ -0.05 & -0.05\end{bmatrix},\ B_2=\begin{bmatrix}-5 & 1\\ -1 & -3\end{bmatrix},\ H_2=\begin{bmatrix}0.1 & 0.12\\ 0.1 & 0.1\end{bmatrix}$$

$$E_2=\begin{bmatrix}0.1 & 0.12\\ 0.12 & 0.1\end{bmatrix},\ E_{d2}=\begin{bmatrix}0.05 & 0.1\\ 0.05 & 0.1\end{bmatrix},\ F_2=\begin{bmatrix}\cos\left(0.5\pi\left(i+j\right)\right) & 0\\ 0 & \cos\left(0.5\pi\left(i+j\right)\right)\end{bmatrix}$$

$$d_h\left(i\right)=2+\sin\left(\frac{\pi i}{2}\right),\ d_v\left(j\right)=2+\sin\left(\frac{\pi j}{2}\right)$$

where state dimensions $n_h=1,\ n_v=1$. From above, it is easy to get that the lower and upper delay bounds along the horizontal and the vertical direction are $d_{hL}=1,\ d_{hH}=3,\ d_{vL}=1,\ d_{vH}=3$. The boundary conditions are as follows

$$x^h\left(i,j\right)=0,\ \ \forall 0\le j\le 20,\ -3\le i<0$$

$$x^h\left(i,j\right)=4,\ \ \forall 0\le j\le 20,\ i=0$$

$$x^h\left(i,j\right)=0,\ \ \forall j>20,\ -3\le i\le 0$$

$$x^v\left(i,j\right)=0,\ \ \forall 0\le i\le 20,\ -3\le j<0$$

$$x^v\left(i,j\right)=3,\ \ \forall 0\le i\le 20,\ j=0$$

$$x^v\left(i,j\right)=0,\ \ \forall i>20,\ -3\le j\le 0$$

The fault matrices $\Omega^k=diag\left\{\omega_{k1},\omega_{k2}\right\},\ k=1,2$

$$\omega_{11}=0.45+0.05\sin\left(0.5\pi\left(i+j\right)\right),\ \omega_{12}=0.55+0.05\sin\left(0.5\pi\left(i+j\right)\right)$$

$$\omega_{21}=0.55+0.05\cos\left(0.5\pi\left(i+j\right)\right),\ \omega_{22}=0.45+0.05\cos\left(0.5\pi\left(i+j\right)\right)$$

Thus $0.4\le\omega_{11}\le 0.5,\ 0.5\le\omega_{12}\le 0.6,\ 0.5\le\omega_{21}\le 0.6,\ 0.4\le\omega_{22}\le 0.5$.

Take $\alpha=0.85$, $\delta_1=\delta_2=0.2$, $\varepsilon_1=\varepsilon_2=0.1$, then solving the matrix inequalities in Theorem 1 gives rise to

$$P^1=\begin{bmatrix}0.5285 & 0\\ 0 & 0.2755\end{bmatrix},\ Q_1^1=\begin{bmatrix}0.0052 & 0\\ 0 & 0.0042\end{bmatrix},\ Q_2^1=\begin{bmatrix}0.0477 & 0\\ 0 & 0.0156\end{bmatrix},$$

$$W_1^1=\begin{bmatrix}0.5470 & 0\\ 0 & 0.3390\end{bmatrix},\ W_2^1=\begin{bmatrix}0.5893 & 0\\ 0 & 0.3581\end{bmatrix},\ \Upsilon^1=\begin{bmatrix}0.0494 & 0.1547\\ 0.1435 & 0.0262\end{bmatrix},$$

$$X_h^1 = \begin{bmatrix} 0.1057 & -0.1047 \\ -0.1047 & 0.1069 \end{bmatrix},\ X_v^1 = \begin{bmatrix} 0.0398 & -0.0390 \\ -0.0390 & 0.0416 \end{bmatrix},\ Y_h^1 = \begin{bmatrix} 0.0137 & -0.0139 \\ -0.0139 & 0.0195 \end{bmatrix},$$

$$Y_v^1 = \begin{bmatrix} 0.0055 & -0.0058 \\ -0.0058 & 0.0084 \end{bmatrix},\ M_h^1 = \begin{bmatrix} -0.1690 \\ 0.1671 \end{bmatrix},\ M_v^1 = \begin{bmatrix} -0.0558 \\ 0.0539 \end{bmatrix},\ N_h^1 = \begin{bmatrix} 0.0001 \\ -0.0352 \end{bmatrix},$$

$$N_v^1 = \begin{bmatrix} 0.0001 \\ -0.0115 \end{bmatrix},\ S_h^1 = \begin{bmatrix} -0.0556 \\ 0.0567 \end{bmatrix},\ S_v^1 = \begin{bmatrix} -0.0180 \\ 0.0188 \end{bmatrix},\ Z = \begin{bmatrix} 0.4141 & 0 \\ 0 & 0.2131 \end{bmatrix}$$

$$P^2 = \begin{bmatrix} 0.4806 & 0 \\ 0 & 0.2875 \end{bmatrix},\ Q_1^2 = \begin{bmatrix} 0.0078 & 0 \\ 0 & 0.0067 \end{bmatrix},\ Q_2^2 = \begin{bmatrix} 0.0426 & 0 \\ 0 & 0.0220 \end{bmatrix},$$

$$W_1^2 = \begin{bmatrix} 0.6218 & 0 \\ 0 & 0.3429 \end{bmatrix},\ W_2^2 = \begin{bmatrix} 0.6634 & 0 \\ 0 & 0.3479 \end{bmatrix},\ X_h^2 = \begin{bmatrix} 0.0759 & -0.0764 \\ -0.0764 & 0.0813 \end{bmatrix},$$

$$X_v^2 = \begin{bmatrix} 0.0555 & -0.0549 \\ -0.0549 & 0.0583 \end{bmatrix},\ Y_h^2 = \begin{bmatrix} 0.0103 & -0.0101 \\ -0.0101 & 0.0155 \end{bmatrix},\ Y_v^2 = \begin{bmatrix} 0.0070 & -0.0074 \\ -0.0074 & 0.0111 \end{bmatrix},$$

$$M_h^2 = \begin{bmatrix} -0.1219 \\ 0.1247 \end{bmatrix},\ M_v^2 = \begin{bmatrix} -0.0627 \\ 0.0609 \end{bmatrix},\ N_h^2 = \begin{bmatrix} 0.0022 \\ -0.0314 \end{bmatrix},\ N_v^2 = \begin{bmatrix} 0.0002 \\ -0.0153 \end{bmatrix},$$

$$S_h^2 = \begin{bmatrix} -0.0395 \\ 0.0431 \end{bmatrix},\ S_v^2 = \begin{bmatrix} -0.0202 \\ 0.0217 \end{bmatrix},\ \Upsilon^2 = \begin{bmatrix} 0.1352 & 0.1098 \\ 0.0851 & -0.0553 \end{bmatrix}.$$

Then $K^1$ and $K^2$ can be obtained by $K^k = \Upsilon^k Z^{-1}$, and the results are as follows

$$K^1 = \begin{bmatrix} 0.1194 & 0.7262 \\ 0.3466 & 0.1228 \end{bmatrix},\ K^2 = \begin{bmatrix} 0.3265 & 0.5155 \\ 0.2055 & -0.2597 \end{bmatrix}$$

Furthermore, we can get $\mu = 3.2292$, From (38), it can be obtained that $\tau_a^* = 7.3$. Choosing $\tau_a = 7.5$, the trajectory of the states $x^h(i, j)$, $x^v(i, j)$ and the corresponding sequence of switching are shown in Figs. 1 and 2.

From Figs. 1 and 2, one can notice that the reliable controller can guarantee the exponential stability and the reliability of the closed-loop system. This demonstrates the effectiveness of the proposed method.

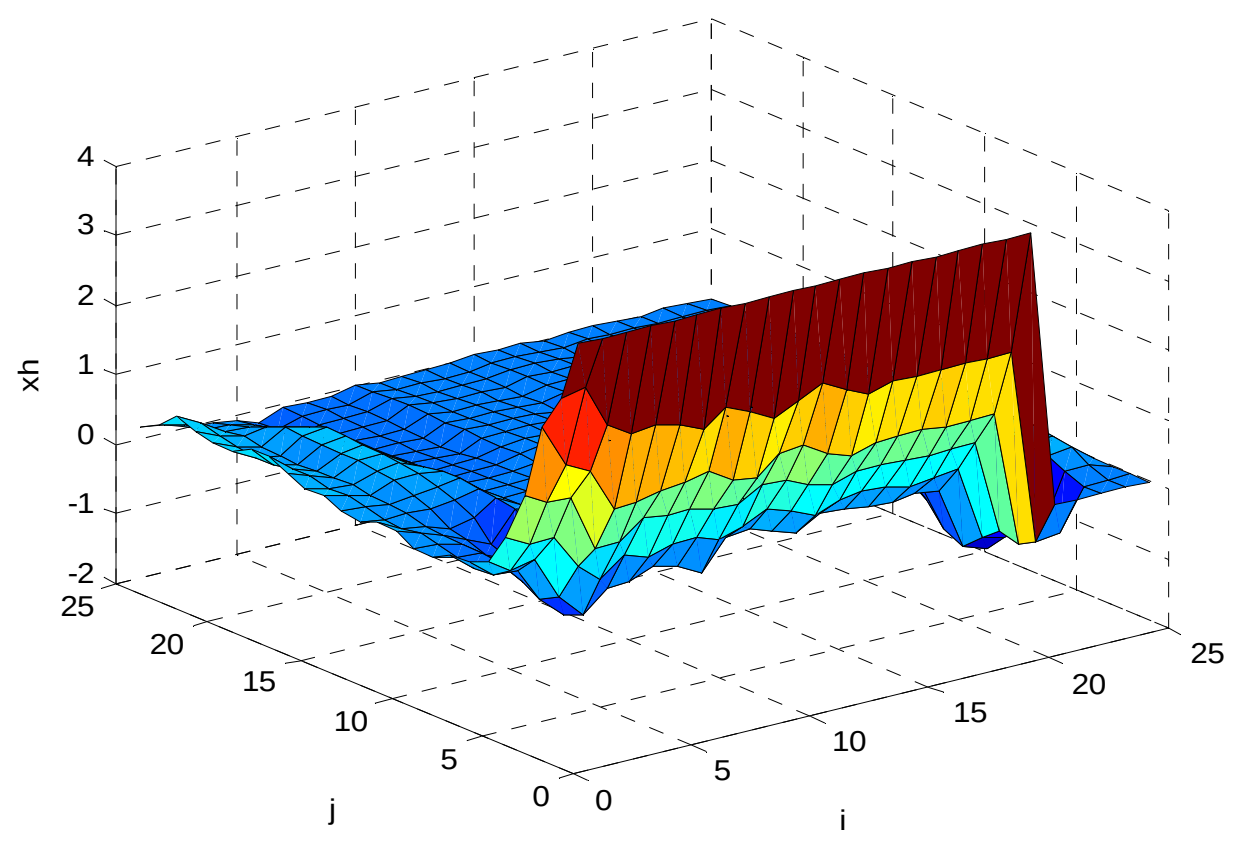

(a) The response of state $x^h(i, j)$

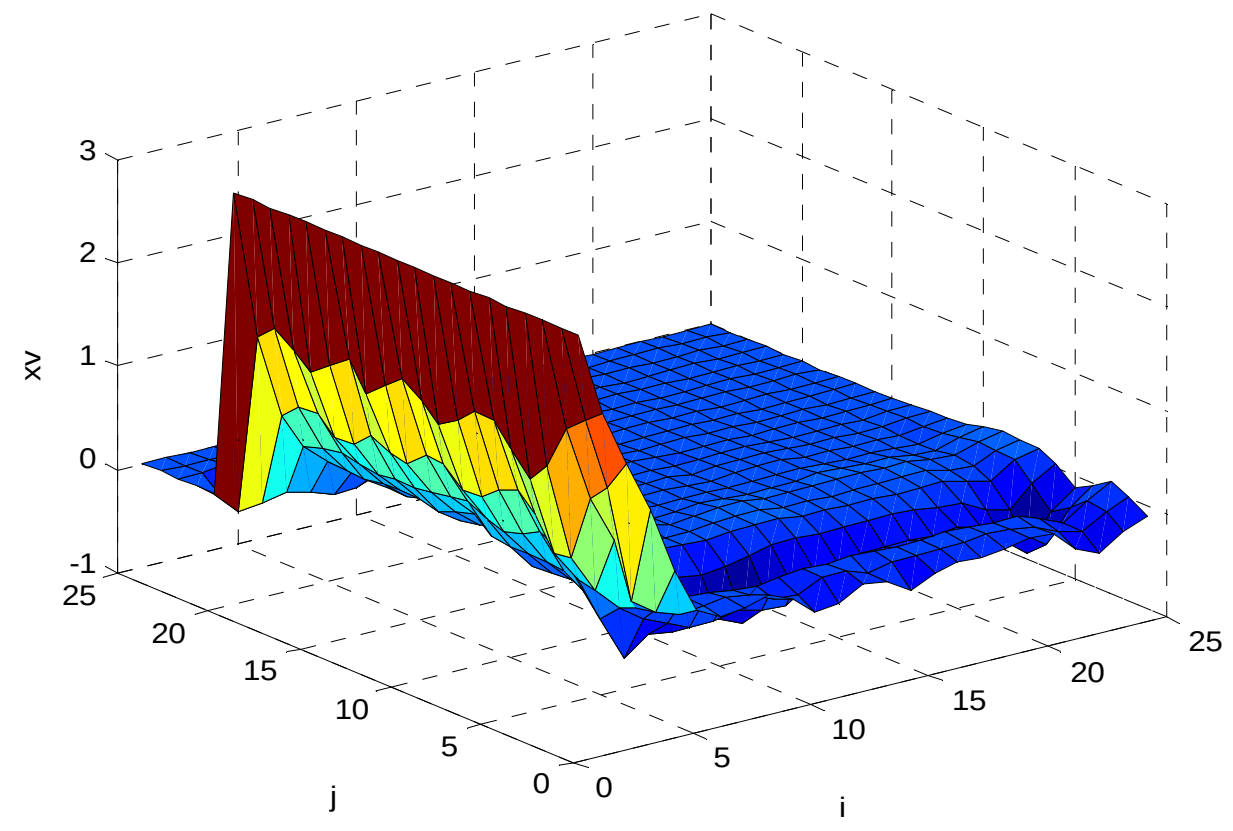

(b) The response of state $x^v(i, j)$

Figure. 1. The trajectory of the states $x^h(i, j)$, $x^v(i, j)$

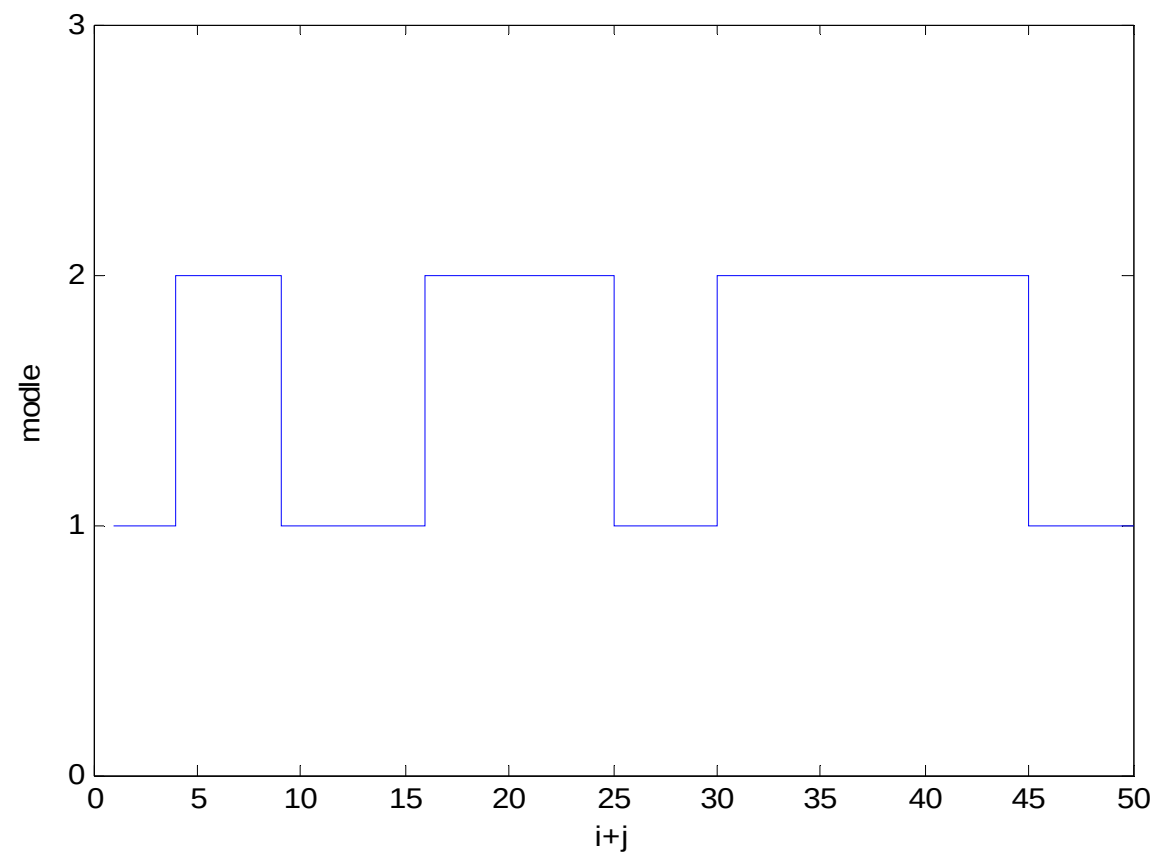

Figure. 2. The corresponding sequence of switching

## 5. Conclusions

This paper has investigated the problem of robust reliable control for a class of uncertain 2D discrete switched systems with actuator failures. A kind of reliable controller design methodology is proposed, and the dwell time approach is utilized for the stability analysis. Sufficient conditions for the existence of such reliable controller are formulated in terms of a set of LMIs. An illustrative example is also given to illustrate the applicability of the proposed approach.

## Acknowledgements

This work was supported by the National Natural Science Foundation of China under Grant No. 60974027 and NUST Research Funding (2011YBXM26).